\documentclass[11pt]{amsart}

\usepackage{amssymb,amsmath,amssymb}
\usepackage{abstract}
\usepackage{amsthm}
\usepackage{mathrsfs}
\usepackage{graphicx}
\usepackage{subfigure}
\usepackage{caption}
\usepackage{float}
\usepackage{epstopdf}
\usepackage{makecell,rotating,multirow,diagbox}
\usepackage{booktabs}
\usepackage{makecell}
\usepackage{color}
\usepackage{url}
\usepackage[colorlinks,linkcolor=blue,anchorcolor=blue,citecolor=red]{hyperref}
\usepackage{appendix}
\usepackage{enumerate}
\usepackage{cite}
\usepackage{graphicx}

\def\d{\mathrm{d}}

\def\bdiv{\mathrm{div}}

\newtheorem{theorem}{Theorem}[section]

\newtheorem{lemma}{Lemma}[section]
\newtheorem{remark}{Remark}[section]
\newtheorem{proposition}{Proposition}[section]
\newtheorem{definition}{Definition}[section]

\begin{document}

	\title[Error Analysis of MIM for High-order Elliptic Equations]{Error Analysis of the Deep Mixed Residual Method for High-order Elliptic Equations}
	\author{Mengjia Bai}
    \address{School of Mathematical Sciences, Soochow University, Suzhou, 215006, China}
    \email{mjbai@stu.suda.edu.cn}
	\author{Jingrun Chen}
    \address{School of Mathematical Sciences and Suzhou Institute for Advanced Research, University of Science and Technology of China, Suzhou, 215123, China}
    \email{jingrunchen@ustc.edu.cn}    
    \author{Rui Du}
    \address{School of Mathematical Sciences and Mathematical Center for Interdisciplinary Research, Soochow University, Suzhou, 215006, China}
    \email{durui@suda.edu.cn}
    \author{Zhiwei Sun}
    \address{Institute of Analysis and Scientific Computing, TU Wien, Wiedner Hauptstra$\ss$e 8–-10, 1040 Wien, Austria}
    \email{zhiwei.sun@asc.tuwien.ac.at}
	\maketitle
	\begin{abstract}
     This paper presents an a priori error analysis of the Deep Mixed Residual method (MIM) for solving high-order elliptic equations with non-homogeneous boundary conditions, including Dirichlet, Neumann, and Robin conditions. We examine MIM with two types of loss functions, referred to as first-order and second-order least squares systems. By providing boundedness and coercivity analysis, we leverage C\'{e}a's Lemma to decompose the total error into the approximation, generalization, and optimization errors. Utilizing the Barron space theory and Rademacher complexity, an a priori error is derived regarding the training samples and network size that are exempt from the curse of dimensionality. Our results reveal that MIM significantly reduces the regularity requirements for activation functions compared to the deep Ritz method, implying the effectiveness of MIM in solving high-order equations.
		 


	\end{abstract}

\section{Introduction}

Partial differential equations (PDEs) are of fundamental importance in modeling phenomena across various disciplines in natural science and society. Developing reliable and efficient numerical methods has a long history in scientific computing and engineering applications. Traditional numerical methods, such as finite difference and finite element, have been successfully established and widely applied. However, these methods often encounter challenges when applied to high-dimensional problems, primarily due to high computational costs. In fact, approximating PDE solutions using traditional methods incurs a computational cost that grows exponentially with the dimensionality of the problem—a phenomenon commonly referred to as the ``curse of dimensionality" (CoD).

In recent years, neural networks have emerged as a promising tool for solving PDEs, demonstrating their potential to address the CoD effectively 
\cite{han2018solving,raissi2019physics,yu2018deep,lyu2020enforcing,lyu2022mim,cai2020deep,zang2020weak}. Notable approaches include the Deep Galerkin method \cite{sirignano2018dgm} and Physics-Informed Neural Networks (PINNs)\cite{raissi2019physics}, which employ the PDE residual in a least-squares framework as the loss function. Another notable approach, the Deep Ritz Method (DRM) \cite{yu2018deep}, leverages the variational form (when available) of the target PDE to define the loss function. More recently, the Deep Mixed Residual method (MIM) \cite{lyu2022mim,lyu2020enforcing} has introduced auxiliary networks to approximate the solution derivatives, allowing for exact enforcement of boundary and initial conditions. Compared to DRM and PINN, MIM has shown advantages in certain models, producing better approximations and accelerating the training process. Additionally, MIM offers unique benefits for handling high-order PDEs by transforming complex high-order problems into lower-order representations, thereby reducing computational complexity and improving solution stability.

In this work, we present an error analysis of the MIM for solving high-order elliptic equations using two-layer neural networks. High-order elliptic equations have extensive applications in materials science \cite{duan2007exact,khain2008generalized}, image processing \cite{andersson1998solution}, and elastic mechanics \cite{landau2012theory}.
To illustrate the MIM framework for high-order equations, consider a $2n$-order elliptic equation with general boundary conditions:
\begin{equation}\label{general model}
    \begin{aligned}
& \Delta^n u=f,&& \quad x \in \Omega, \\
& B(u, \nabla u, \Delta u,\cdots, \nabla\Delta^{n-1} u)= \boldsymbol{g} ,&& \quad x \in \partial \Omega.
\end{aligned}
\end{equation}
MIM introduces auxiliary networks $\phi_i$ and vector-valued networks $\boldsymbol{\psi}_j$ to approximate $\Delta^i u$ and $\nabla \Delta^j u$ for $0\le i,j\le n-1$. Combining the squared residual loss with a penalty term yields the mixed residual loss function:
\begin{equation}\label{first order}
\begin{aligned}
    \|\bdiv\, \boldsymbol{\psi}_{n-1} -f \|_{L^2(\Omega)}^2
    +
    \lambda_1\|B(\phi_0, \boldsymbol{\psi}_0,\cdots, \boldsymbol{\psi}_{n-1}) -\boldsymbol{g} \|_{L^2(\partial\Omega)}^2&\\
    +
    \lambda_2
    \Big(
    \sum_{i=0}^{n-1} 
    \|\nabla \phi_i-\boldsymbol{\psi}_i\|_{L^2(\Omega)}^2
    +
    \sum_{i=0}^{n-2} 
    \| \phi_{i+1} -\bdiv\, \boldsymbol{\psi}_i \|_{L^2(\partial\Omega)}^2
    \Big) &.
    \end{aligned}
\end{equation}
This formulation is also referred to as the first-order least squares system and has been used in the finite element method \cite{cai1994first}. Moreover, we also introduce the second-order least squares system, where we use networks $\varphi_i$ to approximate $\Delta^i u$. Then, the mixed residual loss function is given by:
\begin{equation}\label{second order}
\begin{aligned}
    \|\Delta\, \varphi_{n-1} - f \|_{L^2(\Omega)}^2
    + 
    \lambda_1\|B(\varphi_0, \nabla \varphi_0,\cdots, \nabla\varphi_{n-1}) - \boldsymbol{g} \|_{L^2(\partial\Omega)}^2& \\
    + 
    \lambda_2
    \sum_{i=0}^{n-2} 
    \|\Delta \varphi_i -\varphi_{i+1}\|_{L^2(\Omega)}^2& .
    \end{aligned}
\end{equation}
   
   In this paper, we examine the $2n$-order elliptic equation \eqref{general model} under Dirichlet, Neumann, and Robin boundary conditions. Both first-order and second-order least squares formulations, \eqref{first order}-\eqref{second order}, are analyzed. To derive the error estimates for MIM, we apply a general approach in the notion of bilinear form introduced in \cite{zeinhofer2023unified}, 
   where one can utilize C\'ea's Lemma by performing boundedness and coercivity analysis. As a result, the total error is decomposed into three components: approximation error, generalization error, and optimization error. We bound the approximation error by applying the density of $\mathrm{ReLU}^k$-activated neural networks in Barron space. Moreover, the generalization error is controlled by utilizing Rademacher complexity, while the optimization error is beyond the scope of this article. A significant challenge in our analysis arises from the coercivity analysis of the low-order least squares systems (Equations \eqref{first order} and \eqref{second order}), on how to specifically control the coupling terms. We refer to the work \cite{cai1994first} where the second-order elliptic equation without boundary penalty is studied. In the case of high-order elliptic equations, we employ a perturbation technique, in which we select a special sequence of small parameters to effectively bound these cross terms by decoupled terms. In addition, when considering the Dirichlet boundary condition, it is well-documented \cite{schechter1963p,pmlr-v190-muller22b} that using an $L^2$ boundary penalty leads to a loss of regularity of $3/2$. This implies that approximations in $H^2$ yield a posteriori estimates only in $H^{\frac{1}{2}}$. In this work, we utilize the idea in \cite{li2024priori} and derive a priori estimates in $H^1$ under the Dirichlet boundary condition. This is achieved by extending C\'ea's Lemma and conducting a sup-linear coercivity analysis. 

\subsection{Related works}
Methods including PINN \cite{vahab2022physics}, DRM\cite{yu2018deep}, and MIM \cite{lyu2022mim}  have been utilized to simulate high-order PDEs using neural networks. Additionally, theoretical error analyses have been derived for DRM \cite{siegel2023greedy} and PINN \cite{hong2021priori} when solving high-order elliptic equations, incorporating arbitrary dictionaries of functions and studying greedy algorithms. Furthermore, the works in \cite{zeinhofer2023unified, lu2021priori, mishra2023estimates, mishra2021estimates, muller2021error} provide error analyses for PDEs solved via neural networks. In this study, we concentrate on the error analysis of MIM using two-layer neural networks, leveraging Barron space theory \cite{barron1993universal} to address the curse of dimensionality (CoD). For a deeper understanding of Barron space, we recommend \cite{weinan2019barron, li2024two, lu2021priori}.

\subsection{Our Contributions}
We expand upon the existing literature through the following key contributions: 
\begin{itemize}

\item 
We apply a general approach in the bilinear form introduced in \cite{zeinhofer2023unified} to the MIM, where the error estimate is derived by performing boundedness and coercivity analysis. Our analysis covers first-order and second-order system least squares systems for high-order elliptic equations.

\item 
Our analysis encompasses non-homogeneous boundary conditions, including Dirichlet, Neumann, and Robin conditions, whereas previous works on the DRM \cite{siegel2023greedy} and PINN \cite{hong2021priori} are limited to homogeneous boundary conditions.

\item 
Compared to DRM and PINN, our theoretical analysis reveals that MIM significantly reduces the regularity requirements for activation functions. Specifically, for a $2n$-order elliptic equation, DRM \cite{siegel2023greedy} requires $\mathrm{ReLU}^{n+1}$, while MIM only requires $\mathrm{ReLU}^2$  and $\mathrm{ReLU}^3$ for first-order and second-order least squares systems, respectively.
 
\item 
Previous studies \cite{schechter1963p, pmlr-v190-muller22b} have shown that an $L^2$ Dirichlet boundary penalty results in a regularity loss of $3/2$. In this work, we extend C\'ea’s Lemma by establishing a sup-linear coercivity, yielding an $H^1$ error estimate inspired by the approach in \cite{li2024priori}.

\item 
A core component of our work is the coercivity analysis of low-order systems. The primary challenge lies in how to control the coupling terms, in which we introduce a special sequence of small parameters to effectively control these terms through a perturbation technique.
\end{itemize}

This paper is organized as follows. In Section \ref{Sec: Setup and Main results}, we provide the basic setup and the main results of this work. 
The proof of our main results is given in Section \ref{Sec: The proof of main results} by using a general framework of bilinear form. In Section \ref{Sec:the proof of First-order system}, we present the proof of the crucial coercivity estimates for the first-order least squares system, which is essential for the error estimates. The coercivity estimates for the second-order least squares system are provided in Section
\ref{Sec:the proof of second-order system}.

\section{Main results}\label{Sec: Setup and Main results}
In this section, we introduce the basic setup and present the main results. 

\subsection{Model problem} 
Throughout this paper, we let $\Omega=[0,1]^d$ with dimension $d\geq 1$.
  For any integer $n\geq1$, we consider the following $2n$-order elliptic equation:
\begin{equation}\label{2n-order Laplace equation}
	\Delta ^n u=f\quad	\mbox{in}\ \Omega,
	\end{equation}
with certain boundary conditions, including the Dirichlet boundary:
\begin{equation}\label{Dirichlet boundary condition}
		\big(u,\,\Delta u,\cdots,\,\Delta^{n-1}u\big)
		=
		\boldsymbol{g}_\mathrm{D}\quad \mbox{on}\ \partial\Omega,
	\end{equation}
the Neumann boundary:
\begin{equation}\label{Neumann boundary condition}
 \big(\partial_{\boldsymbol{n}} u,\,
 \partial_{\boldsymbol{n}} \Delta u,\cdots,\,
 \partial_{\boldsymbol{n}} \Delta^{n-1} u\big)
	=
	\boldsymbol{g}_\mathrm{N}\quad \mbox{on}\ \partial\Omega,
\end{equation}
and the Robin boundary:
\begin{equation}\label{Robin boundary condition}
\big(u + \partial_{\boldsymbol{n}} u,\,
 \Delta u + \partial_{\boldsymbol{n}} \Delta u,\cdots,\,
 \Delta^{n-1} u + \partial_{\boldsymbol{n}} \Delta^{n-1} u\big)
=
\boldsymbol{g}_\mathrm{R}\quad \mbox{on}\ \partial\Omega.
\end{equation}
Here, $\boldsymbol{n}$ denotes the outer normal vector of $\partial \Omega$. 

\subsection{Neural networks} 
A neural network is a nonlinear parametric model defined as a concatenation of affine maps and activation functions. 
Let $\boldsymbol{x}\in \mathbb{R}^d$ be the input element and $\boldsymbol{y}\in \mathbb{R}^m$ be the output vector function. An $L$-layer neural network can be written as:
\begin{equation*}
\begin{aligned}
\boldsymbol{h}_0& = \boldsymbol{x},\\
\boldsymbol{h}_\ell
	&=
\sigma(
W_\ell \boldsymbol{h}_{\ell-1} + b_\ell
)
	\quad \mbox{for} \quad
	\ell=1,\cdots, L-1,\\
 \boldsymbol{y}
	&=
W_L \boldsymbol{h}_{L-1} + b_L.
 \end{aligned}
	\end{equation*}
Especially, a neural network is classified as shallow if its depth $L=2$, and as deep otherwise.
Here, $W_\ell$ denotes the weight matrix of the $\ell$-th layer, and $\boldsymbol{b}_\ell$ represents the bias vector. The function $\sigma$ is referred to as an activation function, which introduces non-linearity to enhance the model's expressive power. In this work, we use the activation function $\sigma =\mathrm{ReLU}^k$ defined as follows:
\begin{definition}\label{def: relu}
    Let $k\in \mathbb{N}$, the nonlinear univariate function $\mathrm{ReLU}^k(x)$ is defined by 
    \begin{equation*}
    \mathrm{ReLU}^k(x)
    =
    \left\{
        \begin{aligned}
            &x^k,
            \quad 
            &x\geq 0,\\
            &0,\quad
            &x<0.
        \end{aligned}
        \right.
    \end{equation*}
    In particular, we denote $\mathrm{ReQU}(x) =\mathrm{ReLU}^2(x)$ and $\mathrm{ReCU}(x) =\mathrm{ReLU}^3(x)$. Furthermore, $\mathrm{ReLU}^k$ is applied component-wise, i.e. 
    \begin{equation*}
        \mathrm{ReLU}^k(\boldsymbol{x})
        =
        \big(
        \mathrm{ReLU}^k(x_1),\cdots,
        \mathrm{ReLU}^k(x_n)
        \big),
    \end{equation*}
   for any vector $\boldsymbol{x}=(x_1,\cdots, x_n)$.
\end{definition}

\subsection{Two-layer neural networks and Barron space} 
In this work, we use the two-layer neural networks to approximate Equation \eqref{2n-order Laplace equation}-\eqref{Robin boundary condition}, by applying the loss function in the form of both the first-order least squares system \eqref{first order} and second-order least squares system \eqref{second order}. A suitable functional space as a closed hull of two-layer neural networks is the Barron space. In the following, we introduce the spectral Barron space via a cosine transformation \cite{lu2021priori}. Given a set of cosine functions defined by
\begin{equation}
    \Phi_{\boldsymbol{k}}
    :=
\prod_{i=1}^d \cos(\pi k_i x_i),\quad  \boldsymbol{k} = \{k_i\}_{1\le i\le d},\quad k_i \in \mathbb{N}_0.
\end{equation}
Suppose $u\in L^1(\Omega)$, let $\hat{u}(\boldsymbol{k})$ be the expansion coefficients of $u$ under the basis 
$\left\{\Phi_{\boldsymbol{k}}\right\}_{\boldsymbol{k}\in \mathbb{N}_0^d}$. We define for $s\geq 0$ the spectral Barron space $\mathcal{B}^s(\Omega)$ as
    \begin{equation}\label{def: barron space}
        \mathcal{B}^s(\Omega)
        :=
        \Big\{
       u\in L^1(\Omega): \sum_{\boldsymbol{k}\in \mathbb{N}_0^d}(1+\pi^s|\boldsymbol{k}|_1^s)|\hat{u}(\boldsymbol{k})|
       < \infty
        \Big\}.
    \end{equation}
The spectral Barron norm of a function $u \in\mathcal{B}^s(\Omega)$ is given by
\begin{equation*}
    \|u\|_{\mathcal{B}^s(\Omega)}
    =
    \sum_{\boldsymbol{k}\in \mathbb{N}_0^d}(1+\pi^s|\boldsymbol{k}|_1^s).
\end{equation*}

An important property of Barron functions is that they can be well approximated by two-layer neural networks without the curse of dimensionality, which will be detailed in Lemma \ref{lem: appro error}. To proceed, we assume that the solution of $2n$-order elliptic equation \eqref{2n-order Laplace equation}-\eqref{Robin boundary condition} belongs to the Barron space $\mathcal{B}^{2n+3}$. This can be achieved based on the regularity theory of elliptic equations, assuming that the inhomogeneous term and the boundary term are sufficiently smooth\cite{lions2012non}. In the following, we introduce two models of two-layer neural networks, associated with first-order and second-order least squares systems, respectively. Additionally, their parameters are assumed to be a priori bounded by the Barron norm of the solution. 

\subsubsection{Two-layer neural networks induced by first-order system}
Assume $u^*\in \mathcal{B}^{2n+2}(\Omega)$ is a solution to the $2n$-order elliptic equation \eqref{2n-order Laplace equation}-\eqref{Robin boundary condition}. 
We define a set of {\rm ReQU}-activated networks associated with function $u^*$ by
\begin{equation}\label{def: ReQU activated networks}
\begin{aligned}
    \mathcal{F}_{{\rm ReQU},m}
    :=
    \Big\{
    \boldsymbol{c}&+\sum_{i=1}^m \boldsymbol{a}_i\cdot{\rm ReQU}(W_i\boldsymbol{x}+\boldsymbol{b}_i) \in \mathbb{R}^{n(d+1)} \,\Big|\, \boldsymbol{x}\in \mathbb{R}^d,\\
    & |\boldsymbol{c}|\leq 2\|u^*\|_{\mathcal{B}^{2n+2}},\, 
    |W_i|_2 \leq 1,\,
|\boldsymbol{b}_i|\leq 1,\, 
\sum_{i=1}^m |\boldsymbol{a}_i|\leq 4\|u^*\|_{\mathcal{B}^{2n+2}}   \Big\}.
\end{aligned}
\end{equation}
Denote all the parameters of the network by $\theta$. For any neural network $\boldsymbol{u}_\theta\in  \mathcal{F}_{{\rm  ReQU},m}$, we use the notation $\phi_{i}\in\mathbb{R}$ and $\boldsymbol{\psi}_{i}\in\mathbb{R}^d$ for $i=0,\ldots,n-1$, such that
\begin{equation}\label{representation of neural network}
    \boldsymbol{u}_\theta
    =
    \big(\phi_{0},\,\boldsymbol{\psi}_{0},\cdots,\,\phi_{n-1},\,\boldsymbol{\psi}_{n-1}\big).
\end{equation} 
To simplify the first-order least squares system \eqref{first order}, we introduce some notations. Denote $\boldsymbol{f}=(0,\cdots,0,f)^T\in\mathbb{R}^{n(d+1)}$ and define the matrix operator:
\begin{equation}\label{Def: first-order matrix}
	\mathcal{P}
	: =
	\begin{pmatrix}
		\nabla &                -I_{d\times d} & 0 & \cdots &  0\\
		0      &\mathrm{div}&  -I & \cdots & 0\\
		\vdots& \vdots       &  \vdots& \ddots& \vdots  \\ 
		0      & 0           &         0&     \cdots  &  \mathrm{div}
	\end{pmatrix}_{n(d+1)\times n(d+1)},
\end{equation}
where $I_{d\times d}$ represents $d\times d$ identity matrix.
In addition, given $\boldsymbol{u}_\theta\in  \mathcal{F}_{{\rm  ReQU},m}$, we define the trace operator $S_\alpha$ for $\alpha=\mathrm{D}$, $\mathrm{N}$ and $\mathrm{R}$, representing Dirichlet, Neumann and Robin boundary conditions respectively:
\begin{equation}\label{def: first-order trace}
\left\{
\begin{aligned}
    &S_\mathrm{D} \boldsymbol{u}_\theta
=
\big(\phi_{0},\,\phi_{1},\cdots,\,\phi_{n-1}\big);\\
&S_\mathrm{N} \boldsymbol{u}_\theta
=
\big(\boldsymbol{n}\cdot \boldsymbol{\psi}_{0},\,
\boldsymbol{n}\cdot \boldsymbol{\psi}_{1},\cdots,\,
\boldsymbol{n}\cdot \boldsymbol{\psi}_{n-1}\big);\\
& S_\mathrm{R} \boldsymbol{u}_\theta 
=
\big(S_\mathrm{D} + S_\mathrm{R}\big) \boldsymbol{u}_\theta.
\end{aligned}
\right.
\end{equation}
Then, for Dirichlet and Robin boundary conditions, i.e. $\alpha=\mathrm{D}, \mathrm{R}$, we introduce the expected loss functions regarding the first-order least squares system as below:
\begin{equation}\label{The expected loss functions-first-order system}
\mathcal{L}_\alpha(\boldsymbol{u}_\theta)
=
\|\mathcal{P} \boldsymbol{u}_\theta-\boldsymbol{f}\|_{L^2(\Omega)}^2
+
\lambda\|S_\alpha \boldsymbol{u}_\theta-\boldsymbol{g}_\alpha\|^2_{L^2(\partial\Omega)},
	\end{equation}
 and the corresponding empirical loss function is defined as
\begin{equation}\label{The empirical loss functions-first order system}
\widehat{\mathcal{L}}_\alpha(\boldsymbol{u}_\theta)
=
\frac{|\Omega|}{N}\sum_{j=1}^{N} \big|(\mathcal{P} \boldsymbol{u}_\theta-\boldsymbol{f})(X_j)\big|^2
+
\frac{|\partial\Omega|}{\widehat{N}}\sum_{j=1}^{\widehat{N}}\lambda\big|(S_\alpha \boldsymbol{u}_\theta-\boldsymbol{g}_\alpha)(\widehat{X}_j)\big|^2.
	\end{equation}
    Here, $\{X_j\}_{j=1}^{N} $ and $\{\widehat{X}_j\}_{j=1}^{\widehat{N}} $ are randomly sampled from $\Omega$  and $\partial\Omega$, respectively. 
    The penalty constant $\lambda_2$ in \eqref{first order} is omitted for simplicity, without affecting our main results.
For the Neumann boundary condition, an additional penalty is required in the loss function, since the exact solution of Equation \eqref{2n-order Laplace equation} is unique up to a constant. In this case, we approximate the solution with zero mean, and set the loss function for $\boldsymbol{u}_\theta$ with the representation in \eqref{representation of neural network}:
 \begin{equation}\label{The expected loss functions-first-order system-Neumann}
\mathcal{L}_{\mathrm{N}}(\boldsymbol{u}_\theta)
=
\|\mathcal{P} \boldsymbol{u}_\theta-\boldsymbol{f}\|_{L^2(\Omega)}^2
+
\lambda\|S_{\mathrm{N}} \boldsymbol{u}_\theta-\boldsymbol{g}_{\mathrm{N}}\|^2_{L^2(\partial\Omega)}
+
\mu\big|\int_{\Omega}\phi_{0}\,\mathrm{d}\boldsymbol{x}\big|^2,
	\end{equation}
and the empirical loss function is given by:
 \begin{equation}\label{The empirical loss functions-first order system-Neumann}
 \begin{aligned}
\widehat{\mathcal{L}}_\mathrm{N}(\boldsymbol{u}_\theta)
=&
\frac{|\Omega|}{N}\sum_{j=1}^{N} \big|(\mathcal{P} \boldsymbol{u}_\theta-\boldsymbol{f})(X_j)\big|^2\\
&+
\frac{|\partial\Omega|}{\widehat{N}}\sum_{j=1}^{\widehat{N}}\lambda\big|(S_\mathrm{N} \boldsymbol{u}_\theta-\boldsymbol{g}_\mathrm{N})(\widehat{X}_j)\big|^2
+
\mu\big|\frac{|\Omega|}{N}\sum_{j=1}^N \phi_{0}(X_j)\big|^2.
\end{aligned}
	\end{equation}

\subsubsection{Two-layer neural networks induced by second-order system}
We assume the solutions $u^*\in \mathcal{B}^{2n+3}(\Omega)$, and define a set of {\rm ReCU}-activated networks associated with function $u^*$ as follows:
\begin{equation}\label{def: ReCU activated networks}
\begin{aligned}
    \mathcal{F}_{{\rm ReCU},m}
    :=
    \Big\{
    \boldsymbol{c}
    &+
    \sum_{i=1}^m \boldsymbol{a}_i\cdot{\rm ReCU}(W_i\boldsymbol{x}
    +
    \boldsymbol{b}_i)\in \mathbb{R}^n \, \Big|\, 
    \boldsymbol{x}\in\mathbb{R}^d,
    \\
    &|\boldsymbol{c}|\leq 2\|u^*\|_{\mathcal{B}^{2n+3}},\, |W_i|_2
\leq 1,\, |\boldsymbol{b}_i|\leq 1,\, \sum_{i=1}^m |\boldsymbol{a}_i|\leq 4\|u^*\|_{\mathcal{B}^{2n+3}}   \Big\}.
\end{aligned}
\end{equation}
For any neural network $\boldsymbol{v}_\theta\in  \mathcal{F}_{{\rm  ReCU},m}$, we define the notation $\varphi_{i}\in\mathbb{R}$ for $i=0,\ldots,n-1$, and introduce the following representation:
\begin{equation}\label{representation of v}
    \boldsymbol{v}_\theta
    =
    \big(\varphi_{0},\,\varphi_{1},\cdots,\,\varphi_{n-1}\big).
\end{equation}
To express the loss function in the form of second-order least squares system \eqref{second order}, we denote $\boldsymbol{f}=(0,\cdots,0,f)^T\in\mathbb{R}^{n}$ and define the matrix operator:
\begin{equation}\label{Def: second-order matrix}
	\mathcal{P}^*
	: =
	\begin{pmatrix}
 		\Delta &                -I & 0 & \cdots &  0\\
 		0      &\Delta&  -I & \cdots & 0\\
 		\vdots& \vdots       &  \vdots& \ddots& \vdots  \\ 
 		0      & 0           &         0&     \cdots  &  \Delta
 	\end{pmatrix}_{n\times n}.
\end{equation}
Furthermore, given $\boldsymbol{v}_\theta\in  \mathcal{F}_{{\rm  ReCU},m}$, the trace operator $S_\alpha$ for $\alpha=\mathrm{D}$, $\mathrm{N}$, $\mathrm{R}$, can be defined as follow:
\begin{equation}\label{def: second-order trace}
\left\{
\begin{aligned}
    &S^*_\mathrm{D} \boldsymbol{v}_\theta
=
\big(\varphi_{0},\,\varphi_{1},\cdots,\,\varphi_{n-1}\big);\\
&S^*_\mathrm{N} \boldsymbol{v}_\theta
=
\big(\partial_{\boldsymbol{n}}\varphi_{0},\,
\partial_{\boldsymbol{n}}\varphi_{1},\,
\cdots,\,
\partial_{\boldsymbol{n}}\varphi_{n-1}
\big);\\
& S^*_\mathrm{R} \boldsymbol{v}_\theta 
=
\big(S^*_\mathrm{D} + S^*_\mathrm{R}\big) \boldsymbol{v}_\theta.
\end{aligned}
\right.
\end{equation}
Then, for the case $\alpha=\mathrm{D},\mathrm{R}$, we define the expected loss functions regarding the second-order least squares system \eqref{second order} as
\begin{equation}\label{The expected loss functions-second-order system}
\mathcal{L}_\alpha^*(\boldsymbol{v}_\theta)
=
\|\mathcal{P}^* \boldsymbol{v}_\theta-\boldsymbol{f}\|_{L^2(\Omega)}^2
+
\lambda\|S_\alpha^* \boldsymbol{v}_\theta-\boldsymbol{g}_\alpha\|^2_{L^2(\partial\Omega)},
	\end{equation}
 and the corresponding empirical loss function is given by
\begin{equation}\label{The empirical loss functions-second order system}
\widehat{\mathcal{L}}_\alpha^*(\boldsymbol{v}_\theta)
=
\frac{|\Omega|}{N}\sum_{j=1}^{N}\big|(\mathcal{P}^* \boldsymbol{v}_\theta-\boldsymbol{f})(X_j)\big|^2
+
\frac{|\partial\Omega|}{\widehat{N}}\sum_{j=1}^{\widehat{N}}
\lambda\big|(S_\alpha^* \boldsymbol{v}_\theta-\boldsymbol{g}_\alpha)(\widehat{X}_j)\big|^2.
	\end{equation}
For the Neumann boundary condition, we add a zero-mean penalty and define
 \begin{equation}\label{The expected loss functions-second-order system-Neumann}
\mathcal{L}_{\mathrm{N}}^*(\boldsymbol{v}_\theta)
=
\|\mathcal{P}^* \boldsymbol{v}_\theta-\boldsymbol{f}\|_{L^2(\Omega)}^2
+
\lambda\|S_\mathrm{N}^* \boldsymbol{v}_\theta-\boldsymbol{g}_\alpha\|^2_{L^2(\partial\Omega)}
+
\mu\big| \int_{\Omega}\varphi_{0,\theta}\,\mathrm{d}\boldsymbol{x} \big|^2,
\end{equation}
together with the empirical loss function:
\begin{equation}\label{The empirical loss functions-second order system-Neumann}
\begin{aligned}
\widehat{\mathcal{L}}_\mathrm{N}^*(\boldsymbol{v}_\theta)
=&
\frac{|\Omega|}{N}\sum_{j=1}^{N}\big|(\mathcal{P}^* \boldsymbol{v}_\theta-\boldsymbol{f})(X_j)\big|^2\\
&+
\frac{|\partial\Omega|}{\widehat{N}}\sum_{j=1}^{\widehat{N}}
\lambda\big|(S_\mathrm{N}^* \boldsymbol{v}_\theta-\boldsymbol{g}_\mathrm{N})(\widehat{X}_j)\big|^2
+
\mu\big|\frac{|\Omega|}{N}\sum_{j=1}^N \varphi_{0}(X_j)\big|^2.
\end{aligned}
	\end{equation}

\subsection{Main results}
Now, we present the quantitative estimates for the error between exact solutions and MIM neural networks with Monte Carlo sampling. Define the norm
\begin{equation*}
	\begin{aligned}
		\|\boldsymbol{v}\|_{H(\mathrm{div};\Omega)}
		=
		\big(\|\boldsymbol{v}\|_{L^2(\Omega)}^2
		+
  \|\mathrm{div}\ \boldsymbol{v}\|_{L^2(\Omega)}^2\big)^{\frac{1}{2}}.
	\end{aligned}
\end{equation*}
Moreover, given $N$ and $\widehat{N}$ in \eqref{The empirical loss functions-first order system-Neumann} and \eqref{The empirical loss functions-second order system-Neumann}, which represent the number of sample points in $\Omega$ and on $\partial\Omega$ respectively, we make the assumption that $\widehat{N}=O(\frac{N}{d^2})$.
\begin{theorem}
[First-order least squares system]\label{Thm:main result for first-order hierarch}
    Suppose $u_\alpha^*\in \mathcal{B}^{2n+2}(\Omega)$ with $\alpha=\mathrm{D}, \mathrm{N},\mathrm{R}$ is a solution of problem \eqref{2n-order Laplace equation} satisfying boundary conditions \eqref{Dirichlet boundary condition}, \eqref{Neumann boundary condition} or \eqref{Robin boundary condition}. Let
    \begin{equation*}
        \hat{\boldsymbol{u}}_{\theta}
        =
        \big(\hat{\phi}_0,\,\hat{\boldsymbol{\psi}}_0,\cdots,\hat{\phi}_{n-1},\,\hat{\boldsymbol{\psi}}_{n-1}\big)
        =
        \arg\min_{\boldsymbol{u}_\theta\in \mathcal{F}_{{\rm ReQU},m}}\widehat{\mathcal{L}}_\alpha(\boldsymbol{u}_\theta).
    \end{equation*}
     Then for the cases of Neumann and Robin boundary, i.e. $\alpha=\mathrm{N},\mathrm{R}$, we have
    \begin{equation*}
    \begin{aligned}
        \sum_{k=0}^{n-1}
        \Big(\mathbb{E}\|\hat{\phi}_k -\Delta^{k}u_\alpha^*\|_{H^1(\Omega)}^2
        +
        \mathbb{E}\|\hat{\boldsymbol{\psi}}_{k}&-\nabla(\Delta)^{k}u_\alpha^*\|_{H(\mathrm{div};\Omega)}^2\Big)\\
        \precsim&\,
        \frac{\|u_\alpha^*\|_{\mathcal{B}^{2n+2}(\Omega)}^2}{m}
        +
        \frac{\|u_\alpha^*\|_{\mathcal{B}^{2n+2}(\Omega)}^2}{\sqrt{N}};
        \end{aligned}
    \end{equation*}
    In the case of the Dirichlet boundary, we have
    \begin{equation*}
    \begin{aligned}
        \sum_{k=0}^{n-1}
        \Big(\mathbb{E}\|\hat{\phi}_k-\Delta^{k}u_\mathrm{D}^*\|_{H^1(\Omega)}^2
        +
        \mathbb{E}\|\hat{\boldsymbol{\psi}}_{k}&-\nabla(\Delta)^{k}u_\mathrm{D}^*\|_{H(\mathrm{div};\Omega)}^2\Big)\\
        \precsim&\,
        \frac{\|u_\mathrm{D}^*\|_{\mathcal{B}^{2n+2}(\Omega)}^2}{\sqrt{m}}
        +
        \frac{\|u_\mathrm{D}^*\|_{\mathcal{B}^{2n+2}(\Omega)}^2}{\sqrt[4]{N}},
        \end{aligned}
    \end{equation*}
    where the expectation is taken on the random sampling of training data in $\Omega$ and $\partial\Omega$.
\end{theorem}
\begin{theorem}[Second-order least square system]\label{Thm:main result for second-order hierarch}
    Suppose $u_\alpha^*\in \mathcal{B}^{2n+3}(\Omega)$ with $\alpha=\mathrm{D}, \mathrm{N},\mathrm{R}$ is a solution of problem \eqref{2n-order Laplace equation} satisfying boundary conditions \eqref{Dirichlet boundary condition}, \eqref{Neumann boundary condition} or \eqref{Robin boundary condition}. Let
    \begin{equation*}
        \hat{\boldsymbol{v}}_{\theta}
        =
        \big(\hat{\phi}_0,\cdots,\,\hat{\phi}_{n-1}\big)
        =
        \arg\min_{\boldsymbol{v}_\theta\in \mathcal{F}_{{\rm ReCU},m}}\widehat{\mathcal{L}}_\alpha^*(\boldsymbol{v}_\theta).
    \end{equation*}
    Then for the cases of Neumann and Robin boundary, i.e. $\alpha=\mathrm{N},\mathrm{R}$, we have
    \begin{equation*}
    \begin{aligned}
        \sum_{k=0}^{n-1}\Big(\mathbb{E}\|\hat{\phi}_k-\Delta^{k}u_\alpha^*\|_{H^1(\Omega)}^2
        +
        \mathbb{E}\|\nabla\hat{\phi}_k
        &-\nabla\Delta^{k}u_\alpha^*\|_{H(\mathrm{div};\Omega)}^2
        \Big)\\
        \precsim\,&
        \frac{\|\phi_\alpha^*\|_{\mathcal{B}^{2n+3}(\Omega)}^2}{m}
        +
        \frac{\|\phi_\alpha^*\|_{\mathcal{B}^{2n+3}(\Omega)}^2}{\sqrt{N}};
        \end{aligned}
    \end{equation*}
    In the case of the Dirichlet boundary,  we have
    \begin{equation*}
    \begin{aligned}
       \sum_{k=0}^{n-1}\Big(
       \mathbb{E}\|\hat{\phi}_k-\Delta^{k}u_\mathrm{D}^*\|_{H^1(\Omega)}^2
        +
        \mathbb{E}\|\nabla\hat{\phi}_k
        &-\nabla\Delta^{k}u_\mathrm{D}^*\|_{H(\mathrm{div};\Omega)}^2
        \Big)\\
        \precsim\,&
        \frac{\|u_\mathrm{D}^*\|_{\mathcal{B}^{2n+3}(\Omega)}^2}{\sqrt{m}}
        +
        \frac{\|u_\mathrm{D}^*\|_{\mathcal{B}^{2n+3}(\Omega)}^2}{\sqrt[4]{N}},
        \end{aligned}
    \end{equation*}
    where the expectation is taken on the random sampling of training data over $\Omega$ and $\partial\Omega$.
\end{theorem}

\section{Proof of main results}\label{Sec: The proof of main results}

In this section, we give the detailed proof of Theorem \ref{Thm:main result for first-order hierarch} and Theorem \ref{Thm:main result for second-order hierarch}. The central idea is to use the general approach from \cite{zeinhofer2023unified}, where one can apply C\'ea's Lemma by performing boundedness and coercivity analysis of the bilinear form induced form \eqref{first order} and \eqref{second order}.

\subsection{Coercivity and boundedness}
To apply C\'ea's Lemma, we first introduce the results of coercivity and boundedness. 
\subsubsection{Case I: First-order system} 
Given functions $\boldsymbol{u}, \boldsymbol{w}\in \big(H^1(\Omega)\times H^1(\Omega)^d\big)^n$
\begin{equation}\label{form of U and V}
\begin{aligned}
    \boldsymbol{u} =&(\phi_0,\,\boldsymbol{\psi}_0,\cdots,\,\phi_{n-1},\,\boldsymbol{\psi}_{n-1}),\\
    \boldsymbol{w} =&(\theta_0,\,\boldsymbol{\xi}_0,\cdots,\,\theta_{n-1},\,\boldsymbol{\xi}_{n-1}),
    \end{aligned}
\end{equation}
where $\phi_i,\theta_i \in \mathbb{R}$ and $\boldsymbol{\psi}_i, \boldsymbol{\xi}_i\in \mathbb{R}^d$. Moreover, we introduce the following bilinear form, induced from \eqref{The expected loss functions-first-order system} and \eqref{The expected loss functions-first-order system-Neumann}: for $\alpha=\mathrm{D},\mathrm{R}$, we define
\begin{equation}\label{bilinear form-first order}
	\mathcal{B}_\alpha(\boldsymbol{u},\boldsymbol{w})
	=
	(\mathcal{P}\boldsymbol{u},\mathcal{P}\boldsymbol{w})
	+
 \lambda (S_\alpha \boldsymbol{u},S_\alpha \boldsymbol{w}),
\end{equation}
and for the Neumann boundary condition, we define
\begin{equation}\label{bilinear form-first order-Neumann}
	\mathcal{B}_\mathrm{N}(\boldsymbol{u},\boldsymbol{w})
	=
	(\mathcal{P}\boldsymbol{u},\mathcal{P}\boldsymbol{w})
	+
 \lambda (S_\mathrm{N} \boldsymbol{u},S_\mathrm{N} \boldsymbol{w})
 +
 \mu\int_\Omega \phi_0\,\mathrm{d}\boldsymbol{x}
 \int_\Omega \theta_0\,\mathrm{d}\boldsymbol{x}.
\end{equation}
Here, the matrix operator $\mathcal{P}$ and trace operator $\mathcal{S}_\alpha$ are defined in \eqref{Def: first-order matrix}-\eqref{def: first-order trace}. 

Regarding the coercivity of the bilinear form \eqref{bilinear form-first order}, the case of the second-order elliptic equation with lower-order terms was studied in \cite{li2024priori}. It proved that, for the Neumann boundary, i.e. $\alpha = \mathrm{N}$, the bilinear form \eqref{bilinear form-first order} yields the estimates
\begin{equation*}
    \mathcal{B}_\alpha(\boldsymbol{u},\boldsymbol{u})
    \ge a \| \boldsymbol{u} \|_{H^1}^2,
\end{equation*}
for any function $\boldsymbol{u}\in H^1$. As for the Dirichlet boundary, it follows that
\begin{equation*}
    \mathcal{B}_{\mathrm{D}}(\boldsymbol{u},\boldsymbol{u})
    \ge a \| \boldsymbol{u} \|_{H^{\frac{1}{2}}}^2.
\end{equation*}
Moreover, the $H^{\frac{1}{2}}$ norm on the right-hand side is optimal. In fact, it was shown in \cite{schechter1963p,zeinhofer2023unified} that it is impossible to replace the $H^{\frac{1}{2}}$ by $H^\alpha$ for any $\alpha>\frac{1}{2}$. This implies that the standard analysis cannot provide the coercivity result to derive an $H^1$ error estimate for the Dirichlet boundary. To address this challenge, we follow the idea in \cite{li2024priori}, where a sup-linear coercivity is derived, in the sense that:
\begin{equation*}
    \mathcal{B}_\mathrm{D}(\boldsymbol{u},\boldsymbol{u})
	\ge
 a_E \| \boldsymbol{u} \|_{H^1}^4,
\end{equation*}
for any function $\boldsymbol{u}$ satisfying $\| \boldsymbol{u} \|_{H^1(\Omega)}\le E$, here $E>0$ is to be determined by the exact solution $u^*$, which leads to a priori error estimate. More precisely, we introduce the following result:
\begin{lemma}\label{Lem: coercivity and continuous for first-order}
  For any function $\boldsymbol{u},\,\boldsymbol{w}\in \big(H^1(\Omega)\times H^1(\Omega)^d\big)^n$,  the bilinear form \eqref{bilinear form-first order}, \eqref{bilinear form-first order-Neumann} are bounded for $\alpha=\mathrm{D},\,\mathrm{N},\,\mathrm{R}$, i.e.
\begin{equation*}
		\begin{aligned}
			\mathcal{B}_\alpha(\boldsymbol{u},\boldsymbol{w})
			\le 
			C_B\|\boldsymbol{u}\|_{H^1(\Omega)}
            \|\boldsymbol{w}\|_{H^1(\Omega)}.
		\end{aligned}
	\end{equation*}
 Furthermore, the bilinear form \eqref{bilinear form-first order}, \eqref{bilinear form-first order-Neumann} have coercivity for the Neumann and Robin boundary regarding to $\big(H^1(\Omega)\times H(\mathrm{div};\Omega)\big)^n$, i.e. for $\alpha=\mathrm{N},\,\mathrm{R }$, suppose $\boldsymbol{u}$ satisfies the representation \eqref{form of U and V}, we have
  \begin{equation*}
		\begin{aligned}
			\mathcal{B}_\alpha(\boldsymbol{u},\boldsymbol{u})
			\geq
			C \sum_{k=0}^{n-1}
   \Big( \|\phi_k\|_{H^1(\Omega)}^2 
   +
   \|\boldsymbol{\psi}_{k}\|_{H(\mathrm{div};\Omega)}^2 \Big).
		\end{aligned}
	\end{equation*}
 In addition, if we further assume that $\| \boldsymbol{u} \|_{H^1(\Omega)}\le E$, then it holds the sup-linear coercivity for the Dirichlet boundary, i.e.
 \begin{equation*}
		\begin{aligned}
			\sqrt{\mathcal{B}_\mathrm{D}(\boldsymbol{u},\boldsymbol{u})}
			\geq
			C_E\sum_{k=0}^{n-1}
   \Big( \|\phi_k\|_{H^1(\Omega)}^2 
   +
   \|\boldsymbol{\psi}_{k}\|_{H(\mathrm{div};\Omega)}^2 \Big).
		\end{aligned}
	\end{equation*}
 Here, we use the notation $\|\boldsymbol{u}\|_{H^1(\Omega)}$ to denote the sum of $\|\phi_k\|_{H^1(\Omega)}$ and $\|\boldsymbol{\psi}_{k}\|_{H^1(\Omega)}$.

\end{lemma}
By repeatedly applying the Cauchy-Schwarz and trace inequalities, one can obtain the boundedness of bilinear operator $\mathcal{B}_\alpha$. In addition, the coercivity can be deduced directly by applying Lemma \ref{Thm:Dirichlet boundary condition-first order} together with the boundness results. Lemma \ref{Thm:Dirichlet boundary condition-first order} and Lemma \ref{Thm: second order} are the main ingredients of the analytical part of this work.

\subsubsection{Case II: Second-order system}
For any functions $\boldsymbol{v},\,\boldsymbol{w}\in \big(H^2(\Omega)\big)^n$, we denote
\begin{equation}\label{form of V and W}
\begin{aligned}
    \boldsymbol{v}=&(\varphi_0,\,\varphi_1,\cdots,\,\varphi_{n-1}),\\
    \boldsymbol{w}=&(\rho_0,\,\rho_1,\cdots,\,\rho_{n-1}).
    \end{aligned}
\end{equation}
Define the following bilinear form induced from loss functions \eqref{The expected loss functions-second-order system} and \eqref{The expected loss functions-second-order system-Neumann}: for $\alpha=\mathrm{D},\mathrm{R}$,
\begin{equation}\label{bilinear form-second order}
	\mathcal{B}_\alpha(\boldsymbol{v};\boldsymbol{w})
	=
	(\mathcal{P}^* \boldsymbol{v},\mathcal{P}^* \boldsymbol{w})
	+
     \lambda(S^*_\alpha \boldsymbol{v},S^*_\alpha \boldsymbol{w}),
\end{equation}
and for Neumann boundary condition,
\begin{equation}\label{bilinear form-second order-Neumann}
	\mathcal{B}_\mathrm{N}(\boldsymbol{v};\boldsymbol{w})
	=
	(\mathcal{P}^* \boldsymbol{v},\mathcal{P}^* \boldsymbol{w})
	+
     \lambda(S^*_\mathrm{N} \boldsymbol{v},S^*_\mathrm{N} \boldsymbol{w})
     +
      \mu\int_\Omega \varphi_0\,\mathrm{d}\boldsymbol{x}
 \int_\Omega \rho_0\,\mathrm{d}\boldsymbol{x}.
\end{equation}
Here, 
the matrix operator $\mathcal{P}^*$ and trace operator $\mathcal{S}_\alpha^*$ are defined in \eqref{Def: second-order matrix}-\eqref{def: second-order trace}. 
With the above definitions, we can further obtain results analogous to those for the first-order system.
\begin{lemma}\label{Lem: coercivity and continuous for second-order}
  For any function $\boldsymbol{v},\,\boldsymbol{w}\in \big(H^2(\Omega)\big)^n$,  the bilinear form \eqref{bilinear form-second order}, \eqref{bilinear form-second order-Neumann} are bounded for $\alpha=\mathrm{D},\,\mathrm{N},\,\mathrm{R}$, i.e.
\begin{equation*}
		\begin{aligned}
			\mathcal{B}_\alpha(\boldsymbol{v},\boldsymbol{w})
			\le 
			C_B\|\boldsymbol{v}\|_{H^2(\Omega)}
            \|\boldsymbol{w}\|_{H^2(\Omega)}.
		\end{aligned}
	\end{equation*}
 Moreover, the bilinear form \eqref{bilinear form-second order}, \eqref{bilinear form-second order-Neumann} have coercivity for the Neumann and Robin boundary, i.e. for $\alpha=\mathrm{N},\,\mathrm{R }$, suppose $\boldsymbol{v}$ satisfies representation \eqref{form of V and W}, we have
  \begin{equation*}
		\begin{aligned}
			\mathcal{B}_\alpha(\boldsymbol{v},\boldsymbol{v})
			\geq
			C\sum_{k=0}^{n-1}
   \Big( \|\varphi_k\|_{H^1(\Omega)}^2 
   +
   \|\nabla \varphi_{k}\|_{H(\mathrm{div};\Omega)}^2 \Big). 
		\end{aligned}
	\end{equation*}
 In addition, if we further assume that $\| \boldsymbol{v} \|_{H^2(\Omega)}\le E$, then we have sup-linear coercivity for the Dirichlet boundary, i.e.
 \begin{equation*}
		\begin{aligned}
			\sqrt{\mathcal{B}_\mathrm{D}(\boldsymbol{v},\boldsymbol{v})}
			\geq
			C_E\sum_{k=0}^{n-1}
   \Big( \|\varphi_k\|_{H^1(\Omega)}^2 
   +
   \|\nabla \varphi_{k}\|_{H(\mathrm{div};\Omega)}^2 \Big).
		\end{aligned}
	\end{equation*}
 Here, we use the notation $\|\boldsymbol{v}\|_{H^2(\Omega)}$ to denote the sum of $\|\varphi_k\|_{H^2(\Omega)}$.

\end{lemma}
The boundedness of bilinear operator $\mathcal{B}^*_\alpha$ can be derived directly by using the Cauchy-Schwarz and trace inequalities. The coercivity follows from Lemma \ref{Thm: second order}, together with the boundedness result.

\subsection{Error decomposition}\label{Subsec:generalized abstract frame}
For simplicity, we first introduce an abstract framework proposed in \cite{zeinhofer2023unified}: Given two Hilbert spaces $X,\, Y$, and denote a linear map
\begin{equation*}
	T:X\rightarrow Y,\quad
	\boldsymbol{u}\mapsto T \boldsymbol{u}.
\end{equation*}
Define the corresponding bilinear form for any $\boldsymbol{u},\boldsymbol{v}\in X$:
\begin{equation}\label{general bilinear form}
    \mathcal{B}(\boldsymbol{u},\boldsymbol{v})
	=
	(T\boldsymbol{u},T\boldsymbol{v})_{Y}.
\end{equation}
Furthermore, denote by $\mathcal{F}_\Theta$ the collection
of all neural network functions with a certain parameter space $\Theta$. Then we have the following lemma:
\begin{lemma}[C\'ea's Lemma]\label{Lem: cea lemma}
	For any function $\boldsymbol{f}\in Y$, define the expected loss function $\mathcal{	L}(\boldsymbol{u})=\|T\boldsymbol{u}-\boldsymbol{f}\|_Y^2$, whose unique minimizer is denoted by $\boldsymbol{u}^*$. 
 Assume that the bilinear form $\mathcal{B}$ defined in \eqref{general bilinear form} is bounded, i.e.
 \begin{equation*}
     \mathcal{B}(\boldsymbol{u},\boldsymbol{v})
     \precsim\,
     \|\boldsymbol{u}\|_X\|\boldsymbol{v}\|_X.
 \end{equation*}
 Moreover, there exists a Hilbert space $X\subseteq Z$:
 \begin{itemize}
     \item 
 if $\mathcal{B}$ is coercive with respect to Z, i.e. $\mathcal{B}(\boldsymbol{u},\boldsymbol{u})\ge a \|\boldsymbol{u}\|_Z^2$, 
 then for every $\boldsymbol{u}_\theta\in \mathcal{F}_\Theta$, it holds that 
	\begin{equation}\label{total error estimate}
 \begin{aligned}
		\|\boldsymbol{u}_{\theta}-\boldsymbol{u}^*\|_Z^2
		\precsim\,
		\delta(\boldsymbol{u}_\theta)
		+
		\inf_{\hat{\boldsymbol{u}}_{\theta}\in \mathcal{F}_\Theta}\|\hat{\boldsymbol{u}}_{\theta}-\boldsymbol{u}^*\|_X^2.
  \end{aligned}
	\end{equation}
 \item if $\mathcal{B}$ is sup-linear coercive with respect to Z, i.e. $\mathcal{B}(\boldsymbol{u},\boldsymbol{u})\ge a \|\boldsymbol{u}\|_Z^4$,  then we have
 \begin{equation}\label{total error estimate-super-linear}
 \begin{aligned}
		\|\boldsymbol{u}_{\theta}-\boldsymbol{u}^*\|_Z^2
		\precsim\,
		\sqrt{\delta(\boldsymbol{u}_\theta)}
		+
  \sqrt{
		\inf_{\hat{\boldsymbol{u}}_{\theta}\in \mathcal{F}_\Theta}\|\hat{\boldsymbol{u}}_{\theta}-\boldsymbol{u}^*\|_X^2}.
  \end{aligned}
	\end{equation}
	\end{itemize}
 Here, $	\delta(\boldsymbol{u}_{\theta})=
	\mathcal{L}(\boldsymbol{u}_{\theta})-\inf_{\hat{\boldsymbol{u}}_{\theta}\in \mathcal{F}_\Theta}\mathcal{L}(\hat{\boldsymbol{u}}_{\theta}).$
\end{lemma}

\begin{proof}
	The proof is essentially the same as that in \cite{muller2022error,zeinhofer2023unified}, which is included here for completeness. By the definition \eqref{general bilinear form} of $\mathcal{B}$, we have 
 \begin{equation*}
  \mathcal{	L}(\boldsymbol{u}_\theta)-	\mathcal{	L}(\boldsymbol{u}^*)
  =
     \mathcal{B}(\boldsymbol{u}_\theta-\boldsymbol{u}^*,\boldsymbol{u}_\theta-\boldsymbol{u}^*).
 \end{equation*} 
 Because of the coercivity assumptions of $\mathcal{B}$, there exists a constant $\beta>0$ such that
	\begin{equation*}
		\mathcal{B}(\boldsymbol{u}_\theta-\boldsymbol{u}^*,\boldsymbol{u}_\theta-\boldsymbol{u})
		\geq
		\beta\|\boldsymbol{u}_\theta-\boldsymbol{u}^*\|_Z^k,
	\end{equation*} 
 where $k=2$ or 4 under coercivity or sup-linear coercivity assumptions. 
	On the other hand, the boundedness of $\mathcal{B}$ yields 
	\begin{equation*}
		\begin{aligned}
			\mathcal{	L}(\boldsymbol{u}_\theta)-	\mathcal{	L}(\boldsymbol{u}^*)
			=&
			\mathcal{	L}(\boldsymbol{u}_\theta)
			-
   \inf_{\hat{\boldsymbol{u}}_{\theta}\in \mathcal{F}_\Theta}\mathcal{	L}(\hat{\boldsymbol{u}}_{\theta})
			+
   \inf_{\hat{\boldsymbol{u}}_{\theta}\in \mathcal{F}_\Theta}\big(\mathcal{	L}(\hat{\boldsymbol{u}}_{\theta})-\mathcal{	L}(\boldsymbol{u}^*)\big)\\
			\leq&
			\delta(\boldsymbol{u}_\theta)
   +
   \gamma\inf_{\hat{\boldsymbol{u}}_{\theta}\in \mathcal{F}_\Theta}\|\hat{\boldsymbol{u}}_{\theta}-\boldsymbol{u}^*\|_X^2
		\end{aligned}
	\end{equation*}
 with the constant $\gamma>0$. Hence combining the two estimates and rearranging the terms, we obtain \eqref{total error estimate} and \eqref{total error estimate-super-linear}.
\end{proof}

\begin{remark}\label{remark: decompose}
    Let $\widehat{\mathcal{L}}(\boldsymbol{u})$ be the empirical loss function corrsponding to $\mathcal{L}(\boldsymbol{u})$, we furthermore decompose the error term $\delta(\boldsymbol{u}_\theta)$ as follows:
    \begin{equation*}
    \begin{aligned}
       \delta(\boldsymbol{u}_\theta)
       =&
       \mathcal{L}(\boldsymbol{u}_{\theta})
       -
       \widehat{\mathcal{L}}(\boldsymbol{u}_{\theta})
       +
      \widehat{\mathcal{L}}(\boldsymbol{u}_{\theta})
       -
       \inf_{\hat{\boldsymbol{u}}_{\theta}\in \mathcal{F}_\Theta}
       \widehat{\mathcal{L}}(\hat{\boldsymbol{u}}_{\theta})
       +
       \inf_{\hat{\boldsymbol{u}}_{\theta}\in \mathcal{F}_\Theta}
       \widehat{\mathcal{L}}(\hat{\boldsymbol{u}}_{\theta})
       -
       \inf_{\hat{\boldsymbol{u}}_{\theta}\in \mathcal{F}_\Theta}\mathcal{L}(\hat{\boldsymbol{u}}_{\theta})\\
       \lesssim&
       \sup_{\hat{\boldsymbol{u}}_{\theta}\in \mathcal{F}_\Theta}
       |\mathcal{L}(\hat{\boldsymbol{u}}_{\theta})
       -
       \widehat{\mathcal{L}}(\hat{\boldsymbol{u}}_{\theta})|
       +
       \big(\widehat{\mathcal{L}}(\boldsymbol{u}_{\theta})
       -
       \inf_{\hat{\boldsymbol{u}}_{\theta}\in \mathcal{F}_\Theta}
       \widehat{\mathcal{L}}(\hat{\boldsymbol{u}}_{\theta})\big).
    \end{aligned}
    \end{equation*}
In the last inequality, the first term is called global generalization error, and the second term is called optimization error. Combining with Lemma \ref{Lem: cea lemma}, our estimate is in the form of ``global generalization error+optimization error+approximation error".
\end{remark}

Next, we apply C\'ea's Lemma to the MIM neural networks with first-order and second-order least squares systems.

\subsubsection{Application to first-order system}
In this case, let the spaces $X=\big(H^1(\Omega)\times H^1(\Omega)^d\big)^n$, and we define the linear map $T_{\alpha}$ with $\alpha= \mathrm{D},\mathrm{N},\mathrm{R}$, corresponding to the case of Dirichlet, Neumann and Robin boundary, as
\begin{equation*}
    T_{\alpha} = 
    \left\{\begin{aligned}
        &(\mathcal{P},\, S_{\alpha}), \quad
        \text{when $\alpha= \mathrm{D},\mathrm{R}$},\\
        &(\mathcal{P},\, S_{\alpha},\, \mathcal{I}), \quad
        \text{when $\alpha= \mathrm{N}$},
    \end{aligned}\right.
\end{equation*}
where the matrix operator $\mathcal{P}$ and trace operator $S_{\alpha}$ are defined in \eqref{Def: first-order matrix}, \eqref{def: first-order trace}, and $\mathcal{I}$ denotes the integral on $\Omega$, i.e. $\mathcal{I}u = \int_{\Omega} u\, \d x$. Moreover, we define
\begin{equation}\label{Def: Y space}
    Y_{\alpha} = 
    \left\{\begin{aligned}
        &L^2(\Omega)\oplus L^2(\partial\Omega), \quad
        \text{when $\alpha= \mathrm{D},\mathrm{R}$},\\
        &L^2(\Omega)\oplus L^2(\partial\Omega)\oplus \mathbb{R}, \quad
        \text{when $\alpha= \mathrm{N}$},
    \end{aligned}\right.
\end{equation}
such that $T_{\alpha}:X\rightarrow Y_{\alpha}$. Then the linear map $T_{\alpha}$ induces a bilinear form $\mathcal{B}_{\alpha}$ defined in \eqref{bilinear form-first order}-\eqref{bilinear form-first order-Neumann}, and the corresponding loss function $\mathcal{L}_\alpha$ is given by \eqref{The expected loss functions-first-order system}, \eqref{The expected loss functions-first-order system-Neumann}. Furthermore, set $Z=\big(H^1(\Omega)\times H(\mathrm{div},\Omega)\big)^n$, then Lemma \ref{Lem: coercivity and continuous for first-order} indicates that $\mathcal{B}_{\alpha}$ is coercive with respect to $Z$ when $\alpha= \mathrm{D},\mathrm{R}$, and sup-linear coercive with respect to $Z$ when $\alpha= \mathrm{N}$.

Now we can apply C\'ea's Lemma to the neural networks of first-order least squares system with the spaces and the bilinear form defined above. Suppose $u_\alpha^*\in B^{2n+2}(\Omega)$ is the solution to problem \eqref{2n-order Laplace equation} with three different boundary conditions \eqref{Dirichlet boundary condition}-\eqref{Robin boundary condition}. We use the notation
\begin{equation*}
    \boldsymbol{u}^*_\alpha =
    \big( u^*_\alpha,\, \nabla u^*_\alpha,\cdots,\,\nabla\Delta^{n-1} u^*_\alpha  \big),
\end{equation*}
then it follows that $\boldsymbol{u}^*\in \big(H^1(\Omega)\times H^1(\Omega)^d\big)^n$.
Moreover, Let 
$$\boldsymbol{u}_{\theta} =(\phi_0,\,\boldsymbol{\psi}_0,\cdots,\,\phi_{n-1},\,\boldsymbol{\psi}_{n-1})$$ 
be any neural network in $\mathcal{F}_{{\rm ReQU},m}$.
For the sake of brevity, we denote the approximation error between $\boldsymbol{u}^*$ and $\boldsymbol{u}_{\theta}$, and the generalization error between $\mathcal{L}_{\alpha}$ and $\widehat{\mathcal{L}}_{\alpha}$ by
 \begin{equation}
    \left\{
    \begin{aligned}
        \mathcal{E}_{\rm app}
        =&
        \inf_{\boldsymbol{u}_\theta\in \mathcal{F}_{\mathrm{ReQU},m}}\sum_{k=0}^{n-1}\|\phi_k-\Delta^{k}u_\alpha^*\|_{H^1(\Omega)}^2
        +
        \|\boldsymbol{\psi}_{k}-\nabla(\Delta)^k u_\alpha^*\|_{H^1(\Omega)}^2,\\
       \mathcal{E}_{\rm gen}
       =&
       \sup_{\boldsymbol{u}_\theta\in \mathcal{F}_{\mathrm{ReQU},m}}
       |\mathcal{L}_\alpha(\boldsymbol{u}_\theta)
       -
       \widehat{\mathcal{L}}_{\alpha}(\boldsymbol{u}_\theta)|,
        \end{aligned}
        \right.
    \end{equation}
    where the empirical loss function $\widehat{\mathcal{L}}_{\alpha}$ is given by \eqref{The empirical loss functions-first order system}, \eqref{The empirical loss functions-first order system-Neumann}.
Hence Lemma \ref{Lem: coercivity and continuous for first-order} and C\'ea's Lemma together with Remark \ref{remark: decompose} indicate the result below.

\begin{lemma}[First-order system]\label{Lem:First-order system}
Use the notation defined above, and denote
    \begin{equation*}
        \hat{\boldsymbol{u}}_{\theta}
        =
        \big(\hat{\phi}_0,\,\hat{\boldsymbol{\psi}}_0,\cdots,\hat{\phi}_{n-1},\,\hat{\boldsymbol{\psi}}_{n-1}\big)
        =
        \arg\min_{\boldsymbol{u}_\theta\in \mathcal{F}_{{\rm ReQU},m}}
        \widehat{\mathcal{L}}_{\alpha}(\boldsymbol{u}_\theta).
    \end{equation*}
     For the Neumann and Robin boundary, i.e. $\alpha=\mathrm{N},\mathrm{R}$, we have
    \begin{equation*}
    \begin{aligned}
        \sum_{k=0}^{n-1}\|&\hat{\phi}_k-\Delta^{k}u_\alpha^*\|^2_{H^1(\Omega)}
        +
        \|\hat{\boldsymbol{\psi}}_{k}-\nabla(\Delta)^k u_\alpha^*\|_{H(\mathrm{div};\Omega)}^2
        \lesssim 
       \mathcal{E}_{\rm app}
       +
       \mathcal{E}_{\rm gen}.
        \end{aligned}
    \end{equation*}
     For the Dirichlet boundary, it follows that
    \begin{equation*}
    \begin{aligned}
        \sum_{k=0}^{n-1}\|&\hat{\phi}_k-\Delta^{k}u_\alpha^*\|^2_{H^1(\Omega)}
        +
        \|\hat{\boldsymbol{\psi}}_{k}-\nabla(\Delta)^k u_\alpha^*\|_{H(\mathrm{div};\Omega)}^2
        \lesssim 
        \sqrt{\mathcal{E}_{\rm app}}
        +
         \sqrt{\mathcal{E}_{\rm gen}}.
        \end{aligned}
    \end{equation*}
    
\end{lemma}

\subsubsection{Application to second-order system}
By analogy with the case of the first-order system, we can define the spaces and bilinear form for the second-order system as follows.
Let $X=\big(H^2(\Omega)\big)^n$, and define the linear map $T_\alpha^*$ as
\begin{equation*}
    T_\alpha^* = 
    \left\{\begin{aligned}
        &(\mathcal{P}^*,\, S^*_\alpha), \quad
        \text{when $\alpha=\mathrm{D},\mathrm{R}$},\\
        &(\mathcal{P}^*,\, S^*_\alpha,\, \mathcal{I}), \quad
        \text{when $\alpha=\mathrm{N}$},
    \end{aligned}\right.
\end{equation*}
where the matrix operator $\mathcal{P}^*$ and the trace operator $S^*_\alpha$ are defined in \eqref{Def: second-order matrix} and \eqref{def: second-order trace}.
Therefore  $T^*_\alpha:X\rightarrow Y_\alpha$, where $Y_\alpha$ is the space defined in \eqref{Def: Y space}. 
In addition, the linear map $T^*_\alpha$ induces a bilinear form $\mathcal{B}_{\alpha}$ defined in \eqref{bilinear form-second order}-\eqref{bilinear form-second order-Neumann} with $\alpha= \mathrm{D},\mathrm{N},\mathrm{R}$. 
Define the norm  space $\mathcal{H}$ as follows:
\begin{equation}\label{definition of mathcal H}
	\mathcal{H}(\Omega)=\big\{v\in H^1(\Omega)\,\big|\, \Delta v\in L^2(\Omega)\big\}
\end{equation}
with the norm
\begin{equation*}
	\|v\|_{\mathcal{H}(\Omega)}
	=
	\big(\|v\|_{H^1(\Omega)}^2+\|\Delta v\|_{L^2(\Omega)}^2\big)^{\frac{1}{2}}.
\end{equation*}
 We set $Z=\big(\mathcal{H}(\Omega)\big)^n$, then Lemma \ref{Lem: coercivity and continuous for second-order} indicates that $\mathcal{B}_{\alpha}$ is coercive with respect to $Z$ when $\alpha= \mathrm{D},\mathrm{R}$, and is sup-linear coercive with respect to $Z$ when $\alpha= \mathrm{N}$.

Subsequently, C\'ea's Lemma can be utilized in the neural networks of the second-order least squares system with the spaces and bilinear form defined above. Assume $u_\alpha^*\in \mathcal{B}^{2n+3}(\Omega)$ is the solution to problem \eqref{2n-order Laplace equation} with three different boundary conditions \eqref{Dirichlet boundary condition}-\eqref{Robin boundary condition}. Denote the vector
\begin{equation*}
    \boldsymbol{v}^*_\alpha =
    \big( u^*_\alpha,\,\cdots,\,\Delta^{n-1} u^*_\alpha  \big),
\end{equation*}
then it follows that $\boldsymbol{v}^*_\alpha\in \big(H^2(\Omega)\big)^n$.
 We also define the corresponding approximation error and generalization error in the following manner:
\begin{equation}
    \left\{
    \begin{aligned}
        \mathcal{E}_{\rm app}^*
        =&
         \inf_{\boldsymbol{v}_\theta\in \mathcal{F}_{\mathrm{ReCU},m}}\sum_{k=0}^{n-1}\|\varphi_k-\Delta^{k}u_\alpha^*\|_{H^2(\Omega)}^2;\\
         \mathcal{E}_{\rm gen}^*
         =&
       \sup_{\boldsymbol{v}_\theta\in \mathcal{F}_{\mathrm{ReCU},m}}
       |\mathcal{L}_{\alpha}^*(\boldsymbol{v}_{\theta})
       -
       \widehat{\mathcal{L}}_\alpha^*(\boldsymbol{v}_{\theta})|.
    \end{aligned}
    \right.
\end{equation}

Therefore, by applying Lemma \ref{Lem: coercivity and continuous for first-order} and C$\rm \acute{e}$a Lemma, we can directly infer the following result.
\begin{lemma}[Second-order system]\label{Thm:Second-order system}
Assume  $u_\alpha^*$ is the solution to problem \eqref{2n-order Laplace equation} with three different boundary conditions \eqref{Dirichlet boundary condition}-\eqref{Robin boundary condition}. Let
    \begin{equation*}
        \boldsymbol{v}_{\hat{\theta}}
        =
        \big(\hat{\varphi}_0,\cdots,\,\hat{\varphi}_{n-1}\big)
        =
        \arg\min_{\boldsymbol{v}_\theta\in \mathcal{F}_{{\rm ReCU},m}}\widehat{\mathcal{L}}(\boldsymbol{v}_\theta).
    \end{equation*}
    Then for the Neumann and Robin boundary, i.e. $\alpha=\mathrm{N},\,\mathrm{R}$, one has
    \begin{equation*}
    \begin{aligned}
       \sum_{k=0}^{n-1}\Big(\|\hat{\varphi}_k-\Delta^{k}u_\alpha^*\|_{H^1(\Omega)}^2
        +
        \|\nabla\hat{\varphi}_k
        &-\nabla\Delta^{k}u_\alpha^*\|_{H(\mathrm{div};\Omega)}^2
        \Big)
        \lesssim 
         \mathcal{E}_{\rm app}^*
         +
          \mathcal{E}_{\rm gen}^*,
        \end{aligned}
    \end{equation*}
for the Dirichlet boundary, it follows that
   \begin{equation*}
    \begin{aligned}
        \sum_{k=0}^{n-1}\Big(\|\hat{\varphi}_k-\Delta^{k}u_\mathrm{D}^*\|_{H^1(\Omega)}^2
        +
        \|\nabla\hat{\varphi}_k
        &-\nabla\Delta^{k}u_\mathrm{D}^*\|_{H(\mathrm{div};\Omega)}^2
        \Big)
        \lesssim 
         \sqrt{\mathcal{E}_{\rm app}^*}
         +
          \sqrt{\mathcal{E}_{\rm gen}^*}.
        \end{aligned}
    \end{equation*}
\end{lemma}

From the aforementioned Lemma \ref{Lem:First-order system} and \ref{Thm:Second-order system}, it can be inferred that the error between the neural network approximation and the true solution can be bounded by the approximation error and the generalization error. Therefore, our task below is to estimate these two errors.

\subsection{Approximation error}
Since the neural networks of both first and second-order least squares systems are considered, we introduce the following approximation error regarding ReCU and ReQU-activated neural networks, as defined earlier in \eqref{def: ReQU activated networks} and \eqref{def: ReCU activated networks}. In addition, to address the impact of the curse of dimensionality on approximation errors, our work adopts a method based on range control in Barron space. The approximation property of neural networks using the ReQU activation function has been thoroughly discussed in \cite{li2024priori}, which is presented as follows:
\begin{lemma}\label{lem: appro error}
Given any function $u^*\in \mathcal{B}^{2n+2}(\Omega)$. Let the collection of ReQU-activated neural networks $\mathcal{F}_{{\rm ReQU},m}$ be defined in \eqref{def: ReQU activated networks} associated with $u^*$.
We use the notation $\boldsymbol{u}^* =
    \big( u^*,\, \nabla u^*,\cdots,\,\nabla\Delta^{n-1} u^*  \big)$, then there exists a network $\boldsymbol{u}_m\in \mathcal{F}_{{\rm ReQU},m}$ such that
    \begin{equation*}
        \|\boldsymbol{u}^*-\boldsymbol{u}_m\|_{H^1(\Omega)}^2
        \precsim
        \frac{\|u^*\|_{\mathcal{B}^{2n+2}}^2}{m}.
    \end{equation*}
\end{lemma}
Furthermore, we give the following approximation results regarding ReCU-activated neural networks. Note that the function $u^*$ requires a higher regularity. 
\begin{lemma}\label{Lem: approximation of ReCU network}
        Given any function $u^*\in \mathcal{B}^{2n+3}(\Omega)$. Let the collection of ReCU-activated neural networks $\mathcal{F}_{{\rm ReCU},m}$ be defined in \eqref{def: ReCU activated networks} associated with $u^*$.
We use the notation $\boldsymbol{u}^* =
    \big( u^*,\, \Delta u^*,\cdots,\,\Delta^{n-1} u^*  \big)$,
    then there exists a network $\boldsymbol{u}_m\in \mathcal{F}_{{\rm ReCU},m}$ such that
    \begin{equation*}
        \|\boldsymbol{u}^*-\boldsymbol{u}_m\|_{H^1(\Omega)}^2
        \precsim
        \frac{\|u^*\|_{\mathcal{B}^{2n+3}}^2}{m}.
    \end{equation*}
\end{lemma}
We refer to the work \cite{li2024priori} for the proof of Lemma \ref{lem: appro error}.
The proof of Lemma \ref{Lem: approximation of ReCU network} is shown in Appendix \ref{Sec:Neural-Network Approximation}, where we apply a technique similar to \cite{li2024priori}. More precisely, we construct ReCU-activated networks by utilizing ReLU and ReQU-activated networks. During the construction process, the coefficients rely on the higher regularity of $u^*$, hence in Lemma \ref{Lem: approximation of ReCU network} the assumption $u^*\in \mathcal{B}^{2n+3}(\Omega)$ is necessary. 

A general work on the density of $\mathrm{ReLU}^k$-activated network in Barron space is studied in \cite{li2024two}. Compared to \cite{li2024two}, we provide a sharper bound on the approximation error, where the constants in our estimates are independent of dimension $d$.

\subsection{Generalization error}
 We will illustrate that the quadrature error trained on a finite dataset $\{X_i\}_{i=1}^{N}$ can be estimated by the Rademacher complexity. 
\begin{definition}
 Let $X=\{X_i\}_{i=1}^{N}$ be a set of random variables in $\Omega$ that is independently distributed, and $\varepsilon=\{\varepsilon_i\}_{i=1}^{N}$ be independent Rademacher random variables that take values 
$+1$ or $-1$ with equal probability. Then the \textbf{empirical Rademacher Complexity} of the function class $\mathcal{S}$ is a random variable given by
    \begin{equation*}
        \hat{R}_{N}(\mathcal{S}):=
       \mathbb{E}_\varepsilon\Big[\sup_{f\in \mathcal{S}}\Big|\frac{1}{N}\sum_{i=1}^n\varepsilon_i f(X_i)\Big|\Big].
    \end{equation*}
    Taking its expectation in terms of $X$ yields the \textbf{Rademacher Complexity} of the function class $\mathcal{F}$
    \begin{equation*}
    R_{N}(\mathcal{S})
        =
        \mathbb{E}_X\mathbb{E}_\varepsilon\Big[\sup_{f\in \mathcal{S}}\Big|\frac{1}{N}\sum_{i=1}^{N}\varepsilon_i f(X_i)\Big|\Big].
    \end{equation*}
\end{definition}
In the following, we explore the relation between the generalization error and the Rademacher complexity.
 Given a network $u\in \mathcal{F}$, we denote an expected loss function $\mathcal{L}(u)$ of the form:
\begin{equation}
    \mathcal{L}(u)=\int_\Omega l(u(x))\mathrm{d}\mu(x),
\end{equation}
where $l(y)$ is a function that measures how well the network output $y=u(x)$ fits a given criterion, and $\mu(x)$ is a probability measure.
The following lemma from \cite{li2024priori} fills up the gap between the Rademacher complexity and the quadrature error.
\begin{lemma}\label{Lem: gap between the Rademacher complexity and
the quadrature error}
    Let $\mathcal{F}$ be a set of functions, $X=\{X_i\}_{i=1}^N$ be $i.i.d.$ random variables following the distribution $\mu(x)$. Then
    \begin{equation}
        \mathbb{E}_X\sup_{u\in\mathcal{F}}\Big|\frac{\mathcal{L}(u)}{|\Omega|}-\frac{1}{N}\sum_{i=1}^N l(u(X_i))\Big|
        \leq
        2R_{N}(\mathcal{S}).
    \end{equation}
    Here, $\mathcal{S}:=\big\{l(u)\,|\,u\in\mathcal{F} \big\}$ is referred to as loss function class.
\end{lemma}
Lemma \ref{Lem: gap between the Rademacher complexity and
the quadrature error} indicates that generalization error can be bounded by Rademacher complexity. In the following. we provide an estimation of the Rademacher complexity for the loss function class regarding both first-order and second-order systems.

For the case of the first-order system, let $l(\boldsymbol{u}_\theta)=|\mathcal{P} \boldsymbol{u}_\theta-\boldsymbol{f}|^2$, we define the interior loss function class:
\begin{equation}
    \begin{aligned}
    \mathcal{S}
    =
    \big\{|\mathcal{P} \boldsymbol{u}_\theta-\boldsymbol{f}|^2\,
\big|\, \boldsymbol{u}_\theta\in \mathcal{F}_{{\rm ReQU},m} \big\},
\end{aligned}
\end{equation}
and let $l (\boldsymbol{u}_\theta)=|S_\alpha \boldsymbol{u}_\theta-\boldsymbol{g}_\alpha|^2$. We define the boundary loss function class:
\begin{equation}
    \begin{aligned}
\mathcal{S}_{\alpha}
    =
    \big\{
|S_\alpha \boldsymbol{u}_\theta-\boldsymbol{g}_\alpha|^2\,\big|\, \boldsymbol{u}_\theta\in \mathcal{F}_{{\rm ReQU},m} \big\}.
    \end{aligned}
\end{equation}

By Lemma SM4.7 in \cite{li2024priori}, we can directly get
\begin{lemma}\label{Lem: Rademacher complexity for first-order system} Denote by $N$ and $\widehat{N}$ the number of sample points in $\Omega$ and on $\partial\Omega$, respectively. Assume that $\widehat{N}=O(\frac{N}{d^2})$.
    The function classes $ \mathcal{S},\, \mathcal{S}_{\alpha}$ satisfy
    \begin{equation}
        R_{N}(\mathcal{S})
        +
        R_{\widehat{N}}(\mathcal{S}_{\alpha})
        \leq
        \frac{C\|u^*\|_{\mathcal{B}^{2n+2}}^2}{\sqrt{N}},
    \end{equation}
    where $C$ depends polynomially on dimension $d$.
\end{lemma}

On the other hand, for the case of the second-order system, we similarly define the function classes:
\begin{equation}
    \begin{aligned}
    \mathcal{S}^*
    =
    \big\{|\mathcal{P}^* \boldsymbol{v}_\theta-\boldsymbol{f}|^2\,
\big|\, \boldsymbol{v}_\theta\in \mathcal{F}_{{\rm ReCU},m} \big\},
\end{aligned}
\end{equation}
and boundary loss function class:
\begin{equation}
    \begin{aligned}
\mathcal{S}_{\alpha}^*
    =
    \big\{
|S_\alpha^* \boldsymbol{v}_\theta-\boldsymbol{g}_\alpha|^2\,
\big|\, \boldsymbol{v}_\theta\in \mathcal{F}_{{\rm ReCU},m} \big\}.
    \end{aligned}
\end{equation}

\begin{lemma}\label{Lem: Rademacher complexity for second-order system}
Assume that the number of boundary and interior sample points satisfies $\widehat{N}=O(\frac{N}{d^2})$.
    The function classes $ \mathcal{S}^*,\,\mathcal{S}_{\alpha}^*$ satisfy
    \begin{equation}
        R_{N}(\mathcal{S}^*)
        +
        R_{\widehat{N}}(\mathcal{S}_{\alpha}^*)
        \leq
        \frac{C\|u^*\|_{\mathcal{B}^{2n+3}}^2}{\sqrt{N}},
    \end{equation}
    where $C$ depends polynomially on dimension $d$.
\end{lemma}
To keep the presentation concise, the proof of this lemma will be provided in Appendix \eqref{Subsec: The proof of generalization error}. 

\subsection{Proof of main theorems}
Now, we can utilize the previous results to prove the main theorems. Here we only present the proof of Theorem \ref{Thm:main result for first-order hierarch}, while Theorem \ref{Thm:main result for second-order hierarch} can be proved in the same way.
\begin{proof}[Proof of Theorem \ref{Thm:main result for first-order hierarch}]
   For Neumann and Robin boundary, from Lemma \ref{Lem:First-order system}, we conclude that the total error $\|\hat{\phi}_k -\Delta^{k} u_\alpha^*\|_{H^1(\Omega)}^2$ and $
        \|\hat{\boldsymbol{\psi}}_{k}-\nabla(\Delta)^k u_\alpha^*\|_{H(\mathrm{div};\Omega)}^2$, where $k=0,1,\cdots, n-1$, can be bounded by generalization error
\begin{equation*}
    \mathcal{E}_{\rm gen}
       =
       \sup_{\boldsymbol{u}_\theta\in \mathcal{F}_{\mathrm{ReQU},m}}
       |\mathcal{L}_\alpha(\boldsymbol{u}_\theta)
       -
       \widehat{\mathcal{L}}_{\alpha}(\boldsymbol{u}_\theta)|,
\end{equation*}
        and approximation error
        \begin{equation*}
            \mathcal{E}_{\rm app}
        =
        \inf_{\boldsymbol{u}_\theta\in \mathcal{F}_{\mathrm{ReQU},m}}\sum_{k=0}^{n-1}\|\phi_k-\Delta^{k}u_\alpha^*\|_{H^1(\Omega)}^2
        +
        \|\boldsymbol{\psi}_{k}-\nabla(\Delta)^k u_\alpha^*\|_{H^1(\Omega)}^2.
        \end{equation*}
    By Lemma \ref{lem: appro error}, since $u_\alpha^*\in \mathcal{B}^{2n+2}(\Omega)$, we derive the estimate for the approximation error
    \begin{equation*}
        \mathcal{E}_{\rm app}
        \precsim\,
        \frac{\|u_\alpha^*\|_{\mathcal{B}^{2n+2}(\Omega)}^2}{m}.
    \end{equation*}
Furthermore, for the generalization error, we apply Lemma \ref{Lem: gap between the Rademacher complexity and
the quadrature error} and \ref{Lem: Rademacher complexity for first-order system} to deduce that the estimate holds for $\boldsymbol{u}_\theta\in \mathcal{F}_{\mathrm{ReQU},m}$:
\begin{equation*}
    \mathcal{E}_{\rm gen}
    \precsim \,
    \mathbb{E}|\mathcal{L}_\alpha(\boldsymbol{u}_\theta)
       -
       \widehat{\mathcal{L}}_{\alpha}(\boldsymbol{u}_\theta)|
       \leq
        2|\Omega|\big(R_{N}(\mathcal{S})
        +
        R_{\widehat{N}}(\mathcal{S}_{\alpha})\big)
        \precsim\,
        \frac{\|u^*\|_{\mathcal{B}^{2n+2}}^2}{\sqrt{N}},
\end{equation*}
where the expectation is taken on the random sampling of training data in $\Omega$ and $\partial\Omega$, and we take $\widehat{N}\ge \frac{N}{d^2}$. The constant in the inequality is at most polynomially dependent on $d$. 

Note that for the Neumann case, an additional zero-mean penalty is applied, thus the Rademacher complexity $R_{N}(\mathcal{S}')$ for an extra function class $\mathcal{S}'= \{\boldsymbol{u}_\theta\,
 |\, \boldsymbol{u}_\theta\in \mathcal{F}_{{\rm ReQU},m}  \}$ need to be estimated, which is clearly bounded by  $R_{N}(\mathcal{S})$.

Finally, applying Lemma \ref{Lem:First-order system}, the theorem is proved.  
    
\end{proof}

\section{First-order system for the $2n$-order elliptic equation} \label{Sec:the proof of First-order system}
In this section, we provide the proof of the coercivity estimates for the bilinear form \eqref{bilinear form-first order} and \eqref{bilinear form-first order-Neumann}, induced from the first-order least squares system of the $2n$-order elliptic equation. Suppose $\boldsymbol{u}\in \mathbb{R}^{n(d+1)}$ with the form:
\begin{equation*}
\begin{aligned}
    \boldsymbol{u} =&(\phi_0,\,\boldsymbol{\psi}_0,\cdots,\,\phi_{n-1},\,\boldsymbol{\psi}_{n-1}).
    \end{aligned}
\end{equation*}
\begin{lemma}\label{Thm:Dirichlet boundary condition-first order}
The bilinear form \eqref{bilinear form-first order}-\eqref{bilinear form-first order-Neumann} are coercive with respect to $\big(H^1(\Omega)\times H(\mathrm{div};\Omega)\big)^n$ in the following sense.
For the Neumann and Robin boundary, i.e. $\alpha=\mathrm{N},\mathrm{R}$, there exists a constant $C>0$ such that
	\begin{equation}\label{Neumann boundary condition-first order continuous}
		\begin{aligned}
			\mathcal{B}_\alpha(\boldsymbol{u},\boldsymbol{u})
			\geq
			C\sum_{k=0}^{n-1}
   \Big(\|\boldsymbol{\psi}_{k}\|_{H(\mathrm{div};\Omega)}^2
			+
            \|\phi_k\|_{H^1(\Omega)}^2\Big) ,
		\end{aligned}
	\end{equation}
 and, for the Dirichlet boundary, there exists a constant $C>0$ such that
 \begin{equation}\label{Dirichlet boundary condition-first order continuous}
		\begin{aligned}
		& \Big(	\mathcal{B}^{\frac{1}{2}}_\mathrm{D}(\boldsymbol{u},\boldsymbol{u})
   +
   \sum_{k=0}^{n-1}
            \| \boldsymbol{\psi}_{k} \|_{L^{2}(\partial\Omega)}\Big)
   \mathcal{B}^{\frac{1}{2}}_\mathrm{D}(\boldsymbol{u},\boldsymbol{u})
			\geq
			C\sum_{k=0}^{n-1}
   \Big( \| \boldsymbol{\psi}_{k}\|_{H(\mathrm{div};\Omega)}^2
			+
            \|\phi_k\|_{H^1(\Omega)}^2 \Big).
		\end{aligned}
	\end{equation}
\end{lemma}

\subsection{Some useful inequalies}
In this subsection, we introduce some important inequalities for the proof of our results. Since our work involves various non-homogeneous boundary conditions, let us first introduce the trace inequalities. Functions of $H(\mathrm{div};\Omega)$ admit a well-defined normal trace on $\partial\Omega$. This normal trace $\boldsymbol{v}\cdot\boldsymbol{n}$ lies in $H^{-\frac{1}{2}}(\partial\Omega)$, i.e.
\begin{equation}\label{trace inequality}
	\|\boldsymbol{n}\cdot\boldsymbol{v}\|_{H^{-\frac{1}{2}}(\partial\Omega)}
	\leq 
	\|\boldsymbol{v}\|_{H(\mathrm{div};\Omega)},
	\quad \mbox{for all } \boldsymbol{v}\in H(\mathrm{div};\Omega),
\end{equation}
where $\boldsymbol{n}$ is the outer normal vector. And we have the following inequality from Theorem 1.5 of \cite{girault2012finite},
\begin{equation}\label{trace inequality1}
	\|v\|_{H^{\frac{1}{2}}(\partial\Omega)}
	\leq 
	C_T\|v\|_{H^1(\Omega)},
	\quad \mbox{for all } v\in H^1(\Omega)
\end{equation}
Moreover, we introduce the Poincar\'e-Friedrichs inequality (Theorem 1.9 of \cite{nevcas1967methodes}): 
For any $v\in H^1(\Omega)$, there exists a constant $C_F$ such that 
\begin{equation}\label{Poincaré-Friedrichs inequality}
	\|v\|_{H^1(\Omega)}^2
	\leq 
	C_F\big(\|\nabla v\|_{L^2(\Omega)}^2
	+
 \|v\|_{L^2(\partial\Omega)}^2\big),
\end{equation}
Furthermore, we introduce another Poincar\'e-Friedrichs inequality of the form:
\begin{equation}\label{Poincaré-Friedrichs inequality-2}
	\|v\|_{H^1(\Omega)}^2
	\leq
	C_F\Big(\|\nabla v\|_{L^2(\Omega)}^2
 +
 \big|\int_\Omega v\,\mathrm{d}\boldsymbol{x}\big|^2
 \Big).
\end{equation}

\subsection{Perturbation method}
The difficulty in the coercivity analysis arises from the cross term in $\mathcal{B}_\alpha$. To address this issue, one can use the perturbation method by applying some small coefficients. To this end, we introduce the following parameter series. Given any integer $n\ge 1$, denote $0<\delta<1$, we define the series $0<\delta_k\le 1$, $k=1,2,\dots,2n+1$, by
\begin{equation}\label{def: weighting operator}
    \delta_{k} = \frac{1}{k+1}\delta^{k}.
\end{equation}
Moreover, given $a,b\ge 0$, for $k=2,\dots,2n$, one can apply Young's inequality to derive
\begin{equation*}
\begin{aligned}
    \delta_{k}  ab 
    \le 
    \frac{1}{2} \Big(\delta \sqrt{\frac{k}{k+2}}\Big)\delta_{k}  a^2 
    + 
    \frac{1}{2}\Big(\delta \sqrt{\frac{k}{k+2}}\Big)^{-1} \delta_{k}  b^2 
    = &
    \frac{\epsilon_{k} }{2} \Big( \delta_{k+1}  a^2 
    + 
    \delta_{k-1}   b^2 \Big),
    \end{aligned}
\end{equation*}
where we set 
\begin{equation*}
    \epsilon_{k}  = \Big( \frac{(k+1)^2-1}{(k+1)^2} \Big)^{\frac{1}{2}}.
\end{equation*}
Hence using the monotony of $\epsilon_k$ yields the estimate:
\begin{equation}\label{key estimate of weights}
\begin{aligned}
    \delta_kab 
    \le &
    \frac{\epsilon_{2n}}{2}\Big( \delta_{k+1} a^2 
    + 
    \delta_{k-1}  b^2 \Big).
    \end{aligned}
\end{equation}
Let us define the weighted energy
\begin{equation}\label{weighted energy}
 \begin{aligned}
       \mathcal{E}_{\delta,n}
     =& 
  \sum_{k=0}^{n-1} \delta_{2(n-k)-1}  \| \mathrm{div}\ \boldsymbol{\psi}_{k}\|_{L^2(\Omega)}^2
  + 
  \sum_{k=0}^{n-1}
  \delta_{2(n-k)}  \| \nabla \phi_{k}\|_{L^2(\Omega)}^2,\\
   &+
   \sum_{k=0}^{n-1} 
  \delta_{2(n-k)}  \| \boldsymbol{\psi}_k \|_{L^2(\Omega)}^2 
  +
  \sum_{k=1}^{n-1}\delta_{2(n-k)+1}  \| \phi_{k}\|_{L^2(\Omega)}^2.
  \end{aligned}
\end{equation}
Then we have the following estimate.
\begin{lemma}
For the bilinear form \eqref{bilinear form-first order}-\eqref{bilinear form-first order-Neumann}, we have
     \begin{equation}\label{estimate of B_I n case}
		\begin{aligned}
			\mathcal{B}_\alpha(\boldsymbol{u};\boldsymbol{u})
			\ge&
   (1- \epsilon_{2n} ) \mathcal{E}_{\delta,n}
   - 
   \epsilon_{2n}  \delta_{2n+1}  \| \phi_{0}\|_{L^2(\Omega)}^2
   + \mathcal{T}_{\mathrm{b},\alpha},
		\end{aligned}
	\end{equation}
 where the boundary term $\mathcal{T}_{\mathrm{b},\alpha}$ is defined by
 \begin{equation}\label{boundary term for n case}
 \begin{aligned}
     \mathcal{T}_{\mathrm{b},\alpha}
     =&
     \lambda (S_\alpha \boldsymbol{u},S_\alpha \boldsymbol{u})_{\partial\Omega}
     +
     \tau_\alpha\\
    &-2
  \sum_{k=0}^{n-2}\delta_{2(n-k)-1} (\boldsymbol{n}\cdot \boldsymbol{\psi}_k, \phi_{k+1})_{\partial\Omega}
		-2
  \sum_{k=0}^{n-1}\delta_{2(n-k)} ( \phi_{k},\boldsymbol{n}\cdot \boldsymbol{\psi}_{k})_{\partial\Omega}.
  \end{aligned}
 \end{equation}
 Here, the term $\tau_\alpha$ is the zero-mean penalty that only appears in the Neumann boundary case, i.e.
 \begin{equation*}
 	\tau_\alpha=
 	\left\{\begin{aligned}
 		&\mu\big|\int_\Omega \phi_0\,\mathrm{d}\boldsymbol{x}\big|^2,&& \text{when $\alpha=\mathrm{N}$}, \\
 		& 0, && \text{when $\alpha=\mathrm{D},\mathrm{R}$}.
 	\end{aligned}
 	\right.
 \end{equation*}
\end{lemma}

\begin{proof} 
For convenience, let us denote the matrix operator $\mathcal{P} = \mathcal{P}_1 - \mathcal{P}_2$, where $\mathcal{P}_1$, $\mathcal{P}_2$ are $n(d+1)\times n(d+1)$ matrix given by
\begin{equation*}
    	\mathcal{P}_1
	=
	\begin{pmatrix}
		\nabla & 0 & \cdots &  0 &  0\\
		0      & \bdiv &  \cdots & 0 &  0\\
		\vdots& \vdots       &  \ddots& \vdots & \vdots  \\ 
		0      & 0           &       \cdots  &  \nabla & 0\\
  0      & 0           &       \cdots  & 0 &  \bdiv
	\end{pmatrix},
 \quad
\mathcal{P}_2
	=
	\begin{pmatrix}
		0 &                I_{d\times d} & 0 & \cdots &  0\\
		0      &0&  I & \cdots & 0\\
		\vdots& \vdots       &  \ddots& \ddots& \vdots  \\ 
  0      & 0           &         \cdots&    0  &  I_{d\times d}\\
		0      & 0           &         \cdots&     0  &  0
	\end{pmatrix},
\end{equation*}
Then the bilinear form $\mathcal{B}_\alpha$ in \eqref{bilinear form-first order}-\eqref{bilinear form-first order-Neumann} can be written as
\begin{equation}\label{B_I for n greater than 1}
	\mathcal{B}_\alpha(\boldsymbol{u};\boldsymbol{u})
	=
	\big((\mathcal{P}_1 - \mathcal{P}_2)\boldsymbol{u},\,(\mathcal{P}_1 - \mathcal{P}_2)\boldsymbol{u}\big)
	+
 \lambda(S_\alpha \boldsymbol{u},S_\alpha \boldsymbol{u})_{\partial\Omega}
 +
 \tau_\alpha.
\end{equation} 
Moreover, we apply the parameters \eqref{def: weighting operator} and define the $n(d+1)\times n(d+1)$ weighting matrix:
\begin{equation}
    W_{\delta,n}
	=
	\begin{pmatrix}
		\delta_{2n}I_{d\times d}  & 0 & \cdots &  0 &  0\\
		0      & \delta_{2n-1}&  \cdots & 0 &  0\\
		\vdots& \vdots       &  \ddots& \vdots & \vdots  \\ 
		0      & 0           &       \cdots  &  \delta_{2}I_{d\times d}  & 0\\
  0      & 0           &       \cdots  & 0 &  \delta_{1}
	\end{pmatrix},
\end{equation}
then it follows that
\begin{equation}\label{B_I transform}
    \begin{aligned}
        \mathcal{B}_\alpha(\boldsymbol{u};\boldsymbol{u})
	\ge&
	\big((\mathcal{P}_1 - \mathcal{P}_2)\boldsymbol{u},\,W_{\delta,n} (\mathcal{P}_1 - \mathcal{P}_2)\boldsymbol{u}\big)
	+
\lambda(S_\alpha \boldsymbol{u},S_\alpha \boldsymbol{u})_{\partial\Omega}
+
 \tau_\alpha\\
        = &
        \big(\mathcal{P}_1 \boldsymbol{u},\,W_{\delta,n} \mathcal{P}_1 \boldsymbol{u}\big)
  +
   \big(\mathcal{P}_2 \boldsymbol{u},\,W_{\delta,n} \mathcal{P}_2 \boldsymbol{u}\big)\\
   &-
   2  \big(\mathcal{P}_1 \boldsymbol{u},\,W_\delta \mathcal{P}_2 \boldsymbol{u}\big)
   +
 \lambda(S_\alpha \boldsymbol{u},S_\alpha \boldsymbol{u})_{\partial\Omega}
 +
 \tau_\alpha.
    \end{aligned}
\end{equation}
Here, the first two terms on RHS are positive and are equivalent to the weighted energy we defined in \eqref{weighted energy}:
\begin{equation*}
    \big(\mathcal{P}_1 \boldsymbol{u},\,W_{\delta,n} \mathcal{P}_1 \boldsymbol{u}\big)
  +
   \big(\mathcal{P}_2\boldsymbol{u},\,W_{\delta,n} \mathcal{P}_2 \boldsymbol{u}\big)
   =
    \mathcal{E}_{\delta,n}.
\end{equation*}
Now we estimate the cross term: 
\begin{equation*}
		\begin{aligned}
			-2\big(\mathcal{P}_1 \boldsymbol{u},\,W_{\delta,n} \mathcal{P}_2 \boldsymbol{u}\big)
			=&
             -2\sum_{k=0}^{n-2}\delta_{2(n-k)-1} (\mathrm{div}\ \boldsymbol{\psi}_k, \phi_{k+1})
		-
  2\sum_{k=0}^{n-1}\delta_{2(n-k)} (\nabla \phi_{k},\boldsymbol{\psi}_{k}).
		\end{aligned}
	\end{equation*}
Integrating by parts, one has
 \begin{equation}\label{cross term}
		\begin{aligned}
			&-2\big(\mathcal{P}_1 \boldsymbol{u},\,W_{\delta,n} \mathcal{P}_2 \boldsymbol{u}\big)\\
  =& 2
  \sum_{k=0}^{n-2}\delta_{2(n-k)-1} (\boldsymbol{\psi}_k, \nabla \phi_{k+1})
		+2
  \sum_{k=0}^{n-1}\delta_{2(n-k)} ( \phi_{k},\mathrm{div}\ \boldsymbol{\psi}_{k})\\
  & -2
  \sum_{k=0}^{n-2}\delta_{2(n-k)-1} (\boldsymbol{n}\cdot \boldsymbol{\psi}_k, \phi_{k+1})_{\partial\Omega}
		-2
  \sum_{k=0}^{n-1}\delta_{2(n-k)} ( \phi_{k},\boldsymbol{n}\cdot \boldsymbol{\psi}_{k})_{\partial\Omega}.
		\end{aligned}
	\end{equation}
 Utilizing estimate \eqref{key estimate of weights}, the first two terms can be estimated by
 \begin{equation*}
		\left\{\begin{aligned}
  &2\delta_{2(n-k)-1} (\boldsymbol{\psi}_k, \nabla \phi_{k+1})
  \ge 
  -  \epsilon_{2n}  \Big(
  \delta_{2(n-k)}  \| \boldsymbol{\psi}_k \|_{L^2(\Omega)}^2 
  + 
  \delta_{2(n-k)-2}  \| \nabla \phi_{k+1}\|_{L^2(\Omega)}^2 \Big) ,\\
  &
  2\delta_{2(n-k)} ( \phi_{k},\mathrm{div}\ \boldsymbol{\psi}_{k})
  \ge 
  - \epsilon_{2n}  \Big( \delta_{2(n-k)+1}  \| \phi_{k}\|_{L^2(\Omega)}^2
  +
  \delta_{2(n-k)-1}  \| \mathrm{div}\ \boldsymbol{\psi}_{k}\|_{L^2(\Omega)}^2 \Big).
		\end{aligned}\right.
	\end{equation*}
 Substituting the estimate into \eqref{cross term}, and using notation \eqref{weighted energy}, it follows that
 \begin{equation*}
		\begin{aligned}
			&-2\big(\mathcal{P}_1 \boldsymbol{u},\,W_{\delta,n} \mathcal{P}_2 \boldsymbol{u}\big)\\
			\ge &
             - \epsilon_{2n}  \mathcal{E}_{\delta,n}
   - \epsilon_{2n}  \delta_{2n+1} \| \phi_{0}\|_{L^2(\Omega)}^2\\
  & -2
  \sum_{k=0}^{n-2}\delta_{2(n-k)-1} (\boldsymbol{n}\cdot \boldsymbol{\psi}_k, \phi_{k+1})_{\partial\Omega}
		-2
  \sum_{k=0}^{n-1}\delta_{2(n-k)} ( \phi_{k},\boldsymbol{n}\cdot \boldsymbol{\psi}_{k})_{\partial\Omega}.
		\end{aligned}
	\end{equation*}
 Combining the estimate above and plugging into \eqref{B_I transform} finally leads to \eqref{estimate of B_I n case}.
 \end{proof}
 
\subsection{Estimate of boundary term}
The boundary term $\mathcal{T}_{\mathrm{b},\alpha}$ defined in \eqref{boundary term for n case} can be estimated in the following lemma.
\begin{lemma}
    For the boundary term $\mathcal{T}_{\mathrm{b},\alpha}$, we have
    \begin{itemize}
        \item Neumann boundary. 
        \begin{equation}\label{estimate of boundary Neumann}
            \mathcal{T}_{\mathrm{b},\mathrm{N}}
            \ge -C \delta \mathcal{E}_{\delta,n}
            -
            \epsilon_{2n}  \delta_{2n+1}  \| \phi_{0}\|_{L^2(\Omega)}^2
            +
            \mu\Big|\int_\Omega \phi_0\,\mathrm{d}\boldsymbol{x}\Big|^2.
        \end{equation}
        \item Robin boundary. 
        \begin{equation}\label{estimate of boundary Robin}
            \mathcal{T}_{\mathrm{b},\mathrm{R}}
            \ge \big(2(1-\epsilon_{2n}) - C\delta\big) \sum_{k=0}^{n-1} \delta_{2(n-k)} \Vert \phi_{k}\Vert_{L^2(\partial\Omega)}^2.
        \end{equation}
        \item Dirichlet boundary. 
        \begin{equation}\label{estimate of boundary Dirichlet}
            \mathcal{T}_{\mathrm{b},\mathrm{D}}
            \ge \lambda \sum_{k=0}^{n-1} \Vert \phi_{k}\Vert_{L^2(\partial\Omega)}^2
            -
            \sum_{k=0}^{n-1}
            \| \boldsymbol{\psi}_{k} \|_{L^{2}(\partial\Omega)}
            \| \phi_k\|_{L^{2}(\partial\Omega)},
        \end{equation}
    \end{itemize}
\end{lemma}
\begin{proof}
\noindent\textbf{Case i. Neumann boundary}.
In the case of the Neumann boundary condition,  we apply Young's inequality and estimate the third term in \eqref{boundary term for n case} by
\begin{equation*}
    \begin{aligned}
        -2 \delta_{2(n-k)-1} (\boldsymbol{n}\cdot \boldsymbol{\psi}_k, \phi_{k+1})_{\partial\Omega}
        \ge &
        - 
        \frac{\lambda}{2}
   \|\boldsymbol{n}\cdot\boldsymbol{\psi}_k\|_{L^{2}(\partial\Omega)}^2
			-
            C \delta_{2(n-k)-1}^2\|\phi_{k+1}\|_{L^{2}(\partial\Omega)}^2\\
            \ge &
        - 
        \frac{\lambda}{2}
   \| \boldsymbol{n}\cdot \boldsymbol{\psi}_k\|_{L^{2}(\partial\Omega)}^2
			- 
        C  \delta_{2(n-k)-1}^2 
   \| \phi_{k+1}\|_{H^1(\Omega)}^2,
    \end{aligned}
\end{equation*}
where in the last line, the trace theorem is used. Applying the fact $\delta_{2(n-k)-1}\le \delta_{2(n-k)-2} \le \delta$  for $0\le k\le n-2$, we have
\begin{equation*}
    \begin{aligned}
        -2 &\delta_{2(n-k)-1} (\boldsymbol{n}\cdot \boldsymbol{\psi}_k, \phi_{k+1})_{\partial\Omega}\\
            \ge &
            - 
        \frac{\lambda}{2}
   \| \boldsymbol{n}\cdot \boldsymbol{\psi}_k\|_{L^{2}(\partial\Omega)}^2
        - 
        C \delta \big( \delta_{2(n-k)-2} 
   \| \nabla \phi_{k+1}\|_{L^2(\Omega)}^2
   -
     \delta_{2(n-k)-1} 
   \| \phi_{k+1}\|_{L^2(\Omega)}^2\big).
    \end{aligned}
\end{equation*}
We apply the same argument to the last term in \eqref{boundary term for n case}. For $0\leq k\leq n-1$, we have:
\begin{equation*}
    \begin{aligned}
        -2&
  \delta_{2(n-k)} ( \phi_{k},\boldsymbol{n}\cdot \boldsymbol{\psi}_{k})_{\partial\Omega}\\
        \ge &
        - 
        \frac{\lambda}{2}
   \|\boldsymbol{n}\cdot\boldsymbol{\psi}_{k}\|_{L^{2}(\partial\Omega)}^2
			-
            C \delta_{2(n-k)}^2\|\phi_k\|_{L^{2}(\partial\Omega)}^2\\
            \ge &
        - 
        \frac{\lambda}{2}
   \|\boldsymbol{n}\cdot \boldsymbol{\psi}_{k}\|_{L^{2}(\partial\Omega)}^2
			 - 
        C \delta \big(\delta_{2(n-k)} 
   \| \nabla \phi_k\|_{L^2(\Omega)}^2
   +
   \delta_{2(n-k)+1} 
   \| \phi_k\|_{L^2(\Omega)}^2\big).
    \end{aligned}
\end{equation*}
Here, the inequality $\delta_{2(n-k)}^2\le \delta \delta_{2(n-k)+1}$  is used in the last line. Combining the estimates above, and using notation \eqref{weighted energy}, we obtain
\begin{equation*}
\begin{aligned}
    -2
  \sum_{k=0}^{n-2}\delta_{2(n-k)-1} &(\boldsymbol{n}\cdot \boldsymbol{\psi}_k, \phi_{k+1})_{\partial\Omega}
		-2
  \sum_{k=0}^{n-1}\delta_{2(n-k)} ( \phi_{k},\boldsymbol{n}\cdot \boldsymbol{\psi}_{k})_{\partial\Omega}\\
  & \ge
  -
        \lambda\sum_{k=0}^{n-1}
        \|\boldsymbol{n}\cdot \boldsymbol{\psi}_{k}\|_{L^{2}(\partial\Omega)}^2
        -
        C \delta \mathcal{E}_{\delta,n}
        -
        \epsilon_{2n}  \delta_{2n+1}  \| \phi_{0}\|_{L^2(\Omega)}^2
        .
  \end{aligned}
\end{equation*}
Substituting into \eqref{boundary term for n case}, we can derive the estimate \eqref{estimate of boundary Neumann}.

\noindent\textbf{Case ii. Robin boundary}.
For the case of the Robin boundary condition, we can rewrite the third term in \eqref{boundary term for n case} as
\begin{equation*}
\begin{aligned}
    -2 \delta_{2(n-k)-1} (\boldsymbol{n}\cdot \boldsymbol{\psi}_k, \phi_{k+1})_{\partial\Omega}
        = &
        -2 \delta_{2(n-k)-1} (\boldsymbol{n}\cdot \boldsymbol{\psi}_k + \phi_{k}, \phi_{k+1})_{\partial\Omega}\\
         &+
        2 \delta_{2(n-k)-1} ( \phi_{k}, \phi_{k+1})_{\partial\Omega}.
        \end{aligned}
\end{equation*}
Applying Young's inequality and estimate \eqref{key estimate of weights}, it follows that
\begin{equation*}
    \begin{aligned}
        -2 \delta_{2(n-k)-1} (\boldsymbol{n}\cdot \boldsymbol{\psi}_k, \phi_{k+1})_{\partial\Omega}
        \ge &
        - 
        \frac{\lambda}{2}
   \|\boldsymbol{n}\cdot\boldsymbol{\psi}_k + \phi_{k}\|_{L^{2}(\partial\Omega)}^2
			-
            C \delta_{2(n-k)-1}^2 \|\phi_{k+1}\|_{L^2(\partial\Omega)}^2\\
         &   -
        \epsilon_{2n} \delta_{2(n-k)} \Vert \phi_{k}\Vert_{L^2(\partial\Omega)}^2
        -
        \epsilon_{2n} \delta_{2(n-k)-2} \Vert \phi_{k+1}\Vert_{L^2(\partial\Omega)}^2.
    \end{aligned}
\end{equation*} 
We use the fact $\delta_{2(n-k)-1}\le \delta_{2(n-k)-2} \le \delta$ and sum up the inequality for $k=0,\dots,n-2$, then it follows that
\begin{equation}\label{first term in boundary robin}
    \begin{aligned}
        -2  \sum_{k=0}^{n-2} &\delta_{2(n-k)-1} (\boldsymbol{n}\cdot \boldsymbol{\psi}_k, \phi_{k+1})_{\partial\Omega}\\
        \ge &
        - 
        \frac{\lambda}{2}  \sum_{k=0}^{n-2}
   \|\boldsymbol{n}\cdot\boldsymbol{\psi}_k + \phi_{k}\|_{L^{2}(\partial\Omega)}^2
		-
            (C\delta + 2\epsilon_{2n} )  \sum_{k=0}^{n-1} \delta_{2(n-k)} \|\phi_{k}\|_{L^2(\partial\Omega)}^2.
    \end{aligned}
\end{equation} 
 As for the last term in \eqref{boundary term for n case}, for $0\leq k\leq n-1$, we have:
\begin{equation*}
	\begin{aligned}
    &-2 \delta_{2(n-k)} (\phi_{k}, \boldsymbol{n}\cdot \boldsymbol{\psi}_k)_{\partial\Omega}\\
        = &
        -2 \delta_{2(n-k)} (\phi_{k}, \boldsymbol{n}\cdot \boldsymbol{\psi}_k + \phi_{k} )_{\partial\Omega}
        +
        2 \delta_{2(n-k)} \Vert \phi_{k}\Vert_{L^2(\partial\Omega)}^2\\
        \ge &
        - 
        \frac{\lambda}{2}
        \|\boldsymbol{n}\cdot\boldsymbol{\psi}_k + \phi_{k}\|_{L^{2}(\partial\Omega)}^2
        -
        C \delta_{2(n-k)}^2 \|\phi_k\|_{L^2(\partial\Omega)}^2
        +
        2 \delta_{2(n-k)} \Vert \phi_{k}\Vert_{L^2(\partial\Omega)}^2.
        \end{aligned}
\end{equation*}
Summing up for $k=0,\dots,n-1$, note that $\delta_{2(n-k)} \le \delta$, we can obtain
\begin{equation*}
    \begin{aligned}
        & -
        2 \sum_{k=0}^{n-1}\delta_{2(n-k)} (\phi_{k}, \boldsymbol{n}\cdot \boldsymbol{\psi}_k)_{\partial\Omega}\\
        \ge &
        -  \frac{\lambda}{2} \sum_{k=0}^{n-1}
   \|\boldsymbol{n}\cdot\boldsymbol{\psi}_k + \phi_{k}\|_{L^{2}(\partial\Omega)}^2
            + \sum_{k=0}^{n-1}
        (2 - C\delta) \delta_{2(n-k)} \Vert \phi_{k}\Vert_{L^2(\partial\Omega)}^2.
    \end{aligned}
\end{equation*} 
Combining it with \eqref{boundary term for n case}, \eqref{first term in boundary robin}, we can finally derive \eqref{estimate of boundary Robin}.

\noindent\textbf{Case iii. Dirichlet boundary}.
For the case of the Dirichlet boundary condition, the boundary term in \eqref{boundary term for n case} can be estimated by
\begin{equation*}
\begin{aligned}
    -2 \delta_{2(n-k)-1} (\boldsymbol{n}\cdot \boldsymbol{\psi}_{k+1}, \phi_{k})_{\partial\Omega}
    \ge &
    -2 \delta
   \| \boldsymbol{n}\cdot \boldsymbol{\psi}_{k+1} \|_{L^{2}(\partial\Omega)}
   \| \phi_k\|_{L^{2}(\partial\Omega)},\\
  \delta_{2(n-k)} ( \phi_{k},\boldsymbol{n}\cdot \boldsymbol{\psi}_{k})_{\partial\Omega}
        \ge &
    -2 \delta
   \| \boldsymbol{n}\cdot \boldsymbol{\psi}_{k} \|_{L^{2}(\partial\Omega)}
   \| \phi_k\|_{L^{2}(\partial\Omega)}.
   \end{aligned}
\end{equation*}
Hence the estimate \eqref{estimate of boundary Robin} follows.
\end{proof}

\subsection{Proof of Lemma \ref{Thm:Dirichlet boundary condition-first order}}
Now let us prove the Lemma \ref{Thm:Dirichlet boundary condition-first order} by applying the results above.
For the Neumann case, we substitute \eqref{estimate of boundary Neumann} into \eqref{estimate of B_I n case} to derive that, for $\alpha = \mathrm{N}$,
  \begin{equation}\label{result of Neumann case}
		\begin{aligned}
			\mathcal{B}_\alpha(\boldsymbol{u};\boldsymbol{u})
			\ge&
   (1- \epsilon_{2n} - C\delta ) \mathcal{E}_{\delta,n}
   - 
   2\epsilon_{2n}  \delta_{2n+1}  \| \phi_{0}\|_{L^2(\Omega)}^2
   +
   \mu\Big|\int_\Omega \phi_0\,\mathrm{d}\boldsymbol{x}\Big|^2.
		\end{aligned}
	\end{equation}
We can utilizing the Poincar\'e-Friedrichs inequality \eqref{Poincaré-Friedrichs inequality-2} to deduce
  \begin{equation*}
		\begin{aligned}
			3 \epsilon_{2n}  \delta_{2n+1}  \| \phi_{0}\|_{L^2(\Omega)}^2
			\le &
   3C \epsilon_{2n}  \delta_{2n+1}  \Big(\| \nabla \phi_{0}\|_{L^2(\Omega)}^2
   +
   \Big|\int_\Omega \phi_0\,\mathrm{d}\boldsymbol{x}\Big|^2
   \Big)\\
   = &
   3 C \sqrt{\frac{k-1}{k+1}} \delta  \delta_{2n}  \Big(\| \nabla \phi_{0}\|_{L^2(\Omega)}^2
   +
   \Big|\int_\Omega \phi_0\,\mathrm{d}\boldsymbol{x}\Big|^2
   \Big)\\
   \le &
   C^* \delta \Big(\mathcal{E}_{\delta,n}
   +
   \Big|\int_\Omega \phi_0\,\mathrm{d}\boldsymbol{x}\Big|^2
   \Big).
		\end{aligned}
	\end{equation*}
 Using the result above, it follows from \eqref{result of Neumann case} that
   \begin{equation*}
		\begin{aligned}
			\mathcal{B}_\alpha(\boldsymbol{u};\boldsymbol{u})
			\ge&
   (1- \epsilon_{2n} - C^*\delta ) \mathcal{E}_{\delta,n}
   +
   \epsilon_{2n}  \delta_{2n+1}  \| \phi_{0}\|_{L^2(\Omega)}^2
   +
   (\mu-C^*\delta)\Big|\int_\Omega \phi_0\,\mathrm{d}\boldsymbol{x}\Big|^2.
		\end{aligned}
	\end{equation*}
 Therefore \eqref{Neumann boundary condition-first order continuous} is derived for $\alpha = \mathrm{N}$. 

 For the Robin case, we substitute \eqref{estimate of boundary Robin} into \eqref{estimate of B_I n case} to derive that, for $\alpha = \mathrm{R}$,
  \begin{equation}\label{result of Robin case}
		\begin{aligned}
			\mathcal{B}_\alpha(\boldsymbol{u};\boldsymbol{u})
			\ge&
   (1- \epsilon_{2n} ) \mathcal{E}_{\delta,n}
    - 
   \epsilon_{2n}  \delta_{2n+1}  \| \phi_{0}\|_{L^2(\Omega)}^2\\
   & +
   \big(2(1-\epsilon_{2n}) - C\delta\big) \sum_{k=0}^{n-1} \delta_{2(n-k)} \Vert \phi_{k}\Vert_{L^2(\partial\Omega)}^2.
		\end{aligned}
	\end{equation}
 Note that by Poincar\'e-Friedrichs inequality \eqref{Poincaré-Friedrichs inequality}, we have
 \begin{equation}\label{result of boundary}
 \begin{aligned}
	2 \epsilon_{2n}  \delta_{2n+1} \| \phi_{0}\|_{L^2(\Omega)}^2
	\leq &
 2 C \epsilon_{2n}  \delta_{2n+1}  
	\big(\| \nabla \phi_{0}\|_{L^2(\Omega)}^2
	+
 \|\phi_{0} \|_{L^2(\partial\Omega)}^2\big)\\
	= &
   2 C \sqrt{\frac{k-1}{k+1}} \delta  \delta_{2n}  \big(\| \nabla \phi_{0}\|_{L^2(\Omega)}^2
	+
 \|\phi_{0} \|_{L^2(\partial\Omega)}^2\big)\\
   \le &
   C^* \delta \big( \mathcal{E}_{\delta,n}
	+  \delta_{2n} 
 \|\phi_{0} \|_{L^2(\partial\Omega)}^2\big).
   \end{aligned}
\end{equation}
Substituting into \eqref{result of Robin case}, it leads to
 \begin{equation}
		\begin{aligned}
			\mathcal{B}_\alpha(\boldsymbol{u};\boldsymbol{u})
			\ge&
   (1- \epsilon_{2n} - C^* \delta ) \mathcal{E}_{\delta,n}
    +
   \epsilon_{2n}  \delta_{2n+1}  \| \phi_{0}\|_{L^2(\Omega)}^2\\
   & +
   \big(2(1-\epsilon_{2n}) - C\delta - C^* \delta\big) \sum_{k=0}^{n-1} \delta_{2(n-k)} \Vert \phi_{k}\Vert_{L^2(\partial\Omega)}^2.
		\end{aligned}
	\end{equation}
  Hence \eqref{Neumann boundary condition-first order continuous} is derived for $\alpha = \mathrm{R}$. 

    For the Dirichlet case, one can substitute \eqref{estimate of boundary Dirichlet} into \eqref{estimate of B_I n case} to derive that, for $\alpha = \mathrm{D}$,
  \begin{equation*}
		\begin{aligned}
			\mathcal{B}_\alpha(\boldsymbol{u};\boldsymbol{u})
			\ge&
   (1- \epsilon_{2n} ) \mathcal{E}_{\delta,n}
    - 
   \epsilon_{2n}  \delta_{2n+1}  \| \phi_{0}\|_{L^2(\Omega)}^2\\
   & + \lambda \sum_{k=0}^{n-1} \Vert \phi_{k}\Vert_{L^2(\partial\Omega)}^2
            -
            \sum_{k=0}^{n-1}
            \| \boldsymbol{\psi}_{k} \|_{L^{2}(\partial\Omega)}
            \| \phi_k\|_{L^{2}(\partial\Omega)}.
		\end{aligned}
	\end{equation*}
 Utilizing the result in \eqref{result of boundary}, it follows that
  \begin{equation*}
		\begin{aligned}
			\mathcal{B}_\alpha(\boldsymbol{u};\boldsymbol{u})
			\ge&
   (1- \epsilon_{2n} - C^* \delta ) \mathcal{E}_{\delta,n}
    + 
   \epsilon_{2n}  \delta_{2n+1}  \| \phi_{0}\|_{L^2(\Omega)}^2\\
   & + (\lambda - C^* \delta) \sum_{k=0}^{n-1} \Vert \phi_{k}\Vert_{L^2(\partial\Omega)}^2
            -
            \sum_{k=0}^{n-1}
            \| \boldsymbol{\psi}_{k} \|_{L^{2}(\partial\Omega)}
            \| \phi_k\|_{L^{2}(\partial\Omega)}.
		\end{aligned}
	\end{equation*}
 Moreover, we note the fact
 \begin{equation*}
     \| \phi_k\|_{L^{2}(\partial\Omega)} 
     \le 
     \mathcal{B}^{\frac{1}{2}}_\alpha(\boldsymbol{u};\boldsymbol{u}),
 \end{equation*}
 we finally obtain
  \begin{equation*}
		\begin{aligned}
			&\mathcal{B}_\alpha(\boldsymbol{u};\boldsymbol{u})
   +
   \sum_{k=0}^{n-1}
            \| \boldsymbol{\psi}_{k} \|_{L^{2}(\partial\Omega)}
           \mathcal{B}^{\frac{1}{2}}_\alpha(\boldsymbol{u};\boldsymbol{u})\\
			\ge&
   (1- \epsilon_{2n} - C^* \delta ) \mathcal{E}_{\delta,n}
    + 
   \epsilon_{2n}  \delta_{2n+1}  \| \phi_{0}\|_{L^2(\Omega)}^2
   + (\lambda - C^* \delta) \sum_{k=0}^{n-1} \Vert \phi_{k}\Vert_{L^2(\partial\Omega)}^2.
		\end{aligned}
	\end{equation*}
    
\section{Second-order system for $2n$-order elliptic equation}\label{Sec:the proof of second-order system}
In this section, we give the proof of the coercivity estimates for the bilinear form \eqref{bilinear form-second order} and \eqref{bilinear form-second order-Neumann}, induced from the second-order least squares system of the $2n$-order elliptic equation. Suppose $\boldsymbol{v}\in \mathbb{R}^{n}$ with the form:
\begin{equation*}
\begin{aligned}
    \boldsymbol{v} =&(\varphi_0,\,,\cdots,\,\varphi_{n-1}).
    \end{aligned}
\end{equation*}
\begin{lemma}\label{Thm: second order}
	The bilinear form \eqref{bilinear form-second order}-\eqref{bilinear form-second order-Neumann}  is coercive with respect to $\big(\mathcal{H}(\Omega)\big)^n$ defined in \eqref{definition of mathcal H} in the following sense.
	For the Neumann and Robin boundary, i.e. $\alpha=\mathrm{N},\mathrm{R}$, there exists a constant $C>0$ such that
	\begin{equation}\label{Neumann boundary condition-second order continuous}
		\begin{aligned}
			\mathcal{B}_\alpha(\boldsymbol{v},\boldsymbol{v})
			\geq
			C\sum_{k=0}^{n-1}
   \Big(\|\nabla \varphi_k\|_{H(\mathrm{div};\Omega)}^2
			+
            \|\varphi_k\|_{H^1(\Omega)}^2\Big) .
		\end{aligned}
	\end{equation}
 And for the Dirichlet boundary, there exists a constant $C>0$ such that
 \begin{equation}\label{Dirichlet boundary condition-second order continuous}
		\begin{aligned}
		& \Big(	\mathcal{B}^{\frac{1}{2}}_\mathrm{D}(\boldsymbol{v},\boldsymbol{v})
   +
   \sum_{k=0}^{n-1}
            \| \nabla \varphi_k \|_{L^{2}(\partial\Omega)}\Big)
   \mathcal{B}^{\frac{1}{2}}_\mathrm{D}(\boldsymbol{v},\boldsymbol{v})
			\geq
			C\sum_{k=0}^{n-1}
   \Big( \| \nabla \varphi_k\|_{H(\mathrm{div};\Omega)}^2
			+
            \|\varphi_k\|_{H^1(\Omega)}^2 \Big)
		\end{aligned}
	\end{equation}

\end{lemma}
\subsection{Perturbation method}
Utilizing the notation $\delta_k$ defined in \eqref{def: weighting operator}, we introduce the weighted energy
\begin{equation}\label{define I-second order}
 \begin{aligned}
       \mathcal{E}^*_{\delta,n}
     =& 
  \sum_{k=1}^{n-1}
\delta_{n-k+1}\|\varphi_{k}\|^2_{L^2(\Omega)}
+
\sum_{k=0}^{n-1}
\delta_{n-k}\|\Delta \varphi_{k}\|^2_{L^2(\Omega)}.
  \end{aligned}
\end{equation}
Then we have the following estimate:
\begin{lemma}
For the bilinear form \eqref{bilinear form-second order}, we have
     \begin{equation}\label{estimate of B_I n case second order}
		\begin{aligned}
			\mathcal{B}_\alpha(\boldsymbol{v};\boldsymbol{v})
			\ge&
  (1- \epsilon_{n} ) \mathcal{E}^*_{\delta,n}
   - 
   \epsilon_{n}  \delta_{n+1}  \| \varphi_{0}\|_{L^2(\Omega)}^2
   +
   \mathcal{T}_{\mathrm{b},\alpha},
		\end{aligned}
	\end{equation}
 where the boundary term $\mathcal{T}_{\mathrm{b},\alpha}$ is defined by
 \begin{equation}\label{boundary term for n case second order}
 \begin{aligned}
     \mathcal{T}_{\mathrm{b},\alpha}
     =&
     \lambda(S_\alpha^* \boldsymbol{v},S_\alpha^* \boldsymbol{v})_{\partial\Omega}
     +
     \tau_\alpha\\
    &- 2
  \sum_{k=0}^{n-2}\delta_{n-k} \Big(
   (\partial_{\boldsymbol{n}} \varphi_{k}, \varphi_{k+1})_{\partial\Omega}
		-
   ( \varphi_{k},\partial_{\boldsymbol{n}} \varphi_{k+1})_{\partial\Omega}\Big).
  \end{aligned}
 \end{equation}
  Here, $\tau_\alpha$ is zero-mean penalty that appears in the Neumann boundary case, i.e.
   \begin{equation*}
  	\tau_\alpha=
  	\left\{\begin{aligned}
  		&\mu\big|\int_\Omega \varphi_0\,\mathrm{d}\boldsymbol{x}\big|^2,&& \text{when $\alpha=\mathrm{N}$}, \\
  		& 0, && \text{when $\alpha=\mathrm{D},\mathrm{R}$}.
  	\end{aligned}
  	\right.
  \end{equation*} 
\end{lemma}
\begin{proof}
Let us denote the $n\times n$ matrix 
\begin{equation*}
    	\mathcal{P}_1^*
	=
	\begin{pmatrix}
		\Delta  & 0 & \cdots &  0 &  0\\
		0      & \Delta &  \cdots & 0 &  0\\
		\vdots& \vdots       &  \ddots& \vdots & \vdots  \\ 
		0      & 0           &       \cdots  &  \Delta  & 0\\
  0      & 0           &       \cdots  & 0 &  \Delta
	\end{pmatrix},
 \quad
\mathcal{P}_2^*
	=
	\begin{pmatrix}
		0 &                I & 0 & \cdots &  0\\
		0      &0&  I & \cdots & 0\\
		\vdots& \vdots       &  \ddots& \ddots& \vdots  \\ 
  0      & 0           &         \cdots&    0  &  I\\
		0      & 0           &         \cdots&     0  &  0
	\end{pmatrix}.
\end{equation*}
Then the bilinear form can be written as
\begin{equation}\label{B_I for n greater than 1-second order}
	\mathcal{B}_\alpha(\boldsymbol{v};\boldsymbol{v})
	=
	\big((\mathcal{P}_1^* - \mathcal{P}_2^*)\boldsymbol{v},\,(\mathcal{P}_1^* - \mathcal{P}_2^*)\boldsymbol{v}\big)
	+
 \lambda(S_\alpha^* \boldsymbol{v},S_\alpha^* \boldsymbol{v})_{\partial\Omega}
 +
 \tau_\alpha.
\end{equation} 
Similarly to the proof of the first-order system, we recall parameters $\delta_i$ in \eqref{def: weighting operator} and introduce the  $n\times n$ weighting matrix $W^*_{\delta.n}$ defined by
\begin{equation}
    W^*_{\delta,n}
	=
	\begin{pmatrix}
		\delta_{n}  & 0 & \cdots &  0 &  0\\
		0      & \delta_{n-1}&  \cdots & 0 &  0\\
		\vdots& \vdots       &  \ddots& \vdots & \vdots  \\ 
		0      & 0           &       \cdots  &  \delta_{2}  & 0\\
  0      & 0           &       \cdots  & 0 &  \delta_{1}
	\end{pmatrix}.
\end{equation}
Then one can derives that
\begin{equation}\label{B_I transform-second order}
    \begin{aligned}
        \mathcal{B}_\alpha(\boldsymbol{v};\boldsymbol{v})
	\ge&
	\big((\mathcal{P}_1^* - \mathcal{P}_2^*)\boldsymbol{v},\,W^*_{\delta.n} (\mathcal{P}_1^* - \mathcal{P}_2^*)\boldsymbol{v}\big)
	+
 \lambda(S_\alpha^* \boldsymbol{v},S_\alpha^* \boldsymbol{v})_{\partial\Omega}
 +
 \tau_\alpha\\
        = &
        \big(\mathcal{P}_1^* \boldsymbol{v},\,W^*_{\delta.n} \mathcal{P}_1^* \boldsymbol{v}\big)
  +
   \big(\mathcal{P}_2^*\boldsymbol{v},\,W^*_{\delta.n} \mathcal{P}_2^*\boldsymbol{v}\big)\\
   & -
   2  \big(\mathcal{P}_1^* \boldsymbol{v},\,W^*_{\delta.n} \mathcal{P}_2^*\boldsymbol{v}\big)
   +
 \lambda(S_\alpha^* \boldsymbol{v},S_\alpha^* \boldsymbol{v})_{\partial\Omega}
 +
 \tau_\alpha.
    \end{aligned}
\end{equation}
Here, the first two terms on RHS are equivalent to the weighted energy defined in \eqref{define I-second order}:
\begin{equation*}
    \big(\mathcal{P}_1^* \boldsymbol{v},\,W^*_{\delta.n} \mathcal{P}_1^* \boldsymbol{v}\big)
    +
    \big(\mathcal{P}_2^*\boldsymbol{v},\,W^*_{\delta.n} \mathcal{P}_2^*\boldsymbol{v}\big)
    = 
    \mathcal{E}^*_{\delta,n}.
\end{equation*}
Now we estimate the cross term: 
\begin{equation*}
		\begin{aligned}
			-2
   \big(\mathcal{P}_1^* \boldsymbol{v},\,W^*_{\delta.n} \mathcal{P}_2^*\boldsymbol{v}\big)
			=&
  -2
  \sum_{k=0}^{n-2}\delta_{n-k} (\Delta \varphi_{k}, \varphi_{k+1}).
		\end{aligned}
	\end{equation*}
 Note that by integration by parts, one has
 \begin{equation*}
		\begin{aligned}
			-2\big(\mathcal{P}_1^* \boldsymbol{v},\,W^*_{\delta.n} \mathcal{P}_2^*\boldsymbol{v}\big)
  =& - 2
  \sum_{k=0}^{n-2}\delta_{n-k} \Big( \big(\varphi_{k}, \Delta \varphi_{k+1}\big)	\\
  & +
   \big(\partial_{\boldsymbol{n}} \varphi_{k}, \varphi_{k+1} \big)_{\partial\Omega}
		-
   \big( \varphi_{k},\partial_{\boldsymbol{n}} \varphi_{k+1} \big)_{\partial\Omega}\Big).
		\end{aligned}
	\end{equation*}
 Utilizing estimate \eqref{key estimate of weights}, one can deduce that
 \begin{equation*}
		-2 \delta_{n-k} (\varphi_{k}, \Delta \varphi_{k+1})	
  \geq
  -\epsilon_{n}
  \Big(
\delta_{n-k+1}\|\varphi_{k}\|^2_{L^2(\Omega)}
+
\delta_{n-k-1}\|\Delta \varphi_{k+1}\|^2_{L^2(\Omega)}
  \Big).
\end{equation*}
 Then using the notation \eqref{define I-second order},  it leads to
 \begin{equation*}
		\begin{aligned}
			2\big(\mathcal{P}_1^* \boldsymbol{v},\,W^*_{\delta.n} \mathcal{P}_2^*\boldsymbol{v}\big)
			\ge &
             - \epsilon_{n} \mathcal{E}^*_{\delta,n}
   - 
   \epsilon_{n}  \delta_{n+1} \| \varphi_{0}\|_{L^2(\Omega)}^2\\
  & -2
  \sum_{k=0}^{n-2}\delta_{n-k} \Big(
   (\partial_{\boldsymbol{n}} \varphi_{k}, \varphi_{k+1})_{\partial\Omega}
		-
   ( \varphi_{k},\partial_{\boldsymbol{n}} \varphi_{k+1})_{\partial\Omega}\Big).
		\end{aligned}
	\end{equation*}
 Combining the estimate above, and plugging into \eqref{B_I transform-second order}, it finally leads to
 \begin{equation*}\label{estimate of B_I n case-second order}
		\begin{aligned}
			\mathcal{B}_\alpha(\boldsymbol{v};\boldsymbol{v})
			\ge&
   (1- \epsilon_{n} ) \mathcal{E}^*_{\delta,n}
   - 
   \epsilon_{n}  \delta_{n+1}  \| \varphi_{0}\|_{L^2(\Omega)}^2
   +
    \lambda(S_\alpha^* \boldsymbol{v},S_\alpha^* \boldsymbol{v})_{\partial\Omega}
    +
    \tau_\alpha\\
    &- 2
  \sum_{k=0}^{n-2}\delta_{n-k} \Big(
   (\partial_{\boldsymbol{n}} \varphi_{k}, \varphi_{k+1})_{\partial\Omega}
		-
   ( \varphi_{k},\partial_{\boldsymbol{n}} \varphi_{k+1})_{\partial\Omega}\Big).
		\end{aligned}
	\end{equation*}
 The Lemma is proved.
 \end{proof}

\subsection{Estimate of boundary term}
Now let us estimate the boundary term $\mathcal{T}_{\mathrm{b},\alpha}$ in \eqref{boundary term for n case second order}.
\begin{lemma}
    For the boundary term $\mathcal{T}_{\mathrm{b},\alpha}$, we have
    \begin{itemize}
        \item Neumann boundary. 
        \begin{equation}\label{boundary estimate of D case n greater than 1-second order}
\begin{aligned}
     \mathcal{T}_{\mathrm{b},\mathrm{N}}
  \ge& 
  - 
        C\delta 
       \mathcal{E}^*_{\delta,n}
       -
       \epsilon_{n}  \delta_{n+1}  \| \varphi_{0}\|_{L^2(\Omega)}^2\\
        &+
        \Big(\frac{\lambda}{2} - C \delta \Big)
        \sum_{k=0}^{n-1}
        \|\partial_{\boldsymbol{n}} \varphi_{k}\|_{L^{2}(\partial\Omega)}^2
        +
        \mu\Big|\int_\Omega \varphi_0\,\mathrm{d}\boldsymbol{x}\Big|^2.
  \end{aligned}
\end{equation}
        \item Robin boundary. 
        \begin{equation}\label{estimate of boundary Robin second order}
            \mathcal{T}_{\mathrm{b},\mathrm{R}}
  \ge
  - 
        C \delta 
       \mathcal{E}^*_{\delta,n}
       +
      \Big(\frac{\lambda}{2} -  C \delta \Big)
       \sum_{k=0}^{n-1}
        \|\varphi_{k} + \partial_{\boldsymbol{n}} \varphi_k\|_{L^{2}(\partial\Omega)}^2.
        \end{equation}
        \item Dirichlet boundary. 
        \begin{equation}\label{estimate of boundary Dirichlet second order}
            \mathcal{T}_{\mathrm{b},\mathrm{D}}
            \ge \lambda \sum_{k=0}^{n-1} \Vert \varphi_{k}\Vert_{L^2(\partial\Omega)}^2
            -
            \sum_{k=0}^{n-1}
            \| \partial_{\boldsymbol{n}} \varphi_{k+1} \|_{L^{2}(\partial\Omega)}
            \| \varphi_k\|_{L^{2}(\partial\Omega)}.
        \end{equation}
    \end{itemize}
\end{lemma}
\begin{proof}
\noindent\textbf{Case i. Neumann boundary}.
For the Neumann boundary case, we utilize Young's inequality and trace theorem for the third term in \eqref{boundary term for n case second order} and derive the estimate
\begin{equation*}
    \begin{aligned}
        2 \delta_{n-k} (\partial_{\boldsymbol{n}} \varphi_{k}, \varphi_{k+1})_{\partial\Omega}
        \ge &
        - 
        C \delta_{n-k}^2
   \|\varphi_{k+1}\|_{L^{2}(\partial\Omega)}^2
			-
            \frac{\lambda}{4}
            \|\partial_{\boldsymbol{n}} \varphi_{k}\|_{L^{2}(\partial\Omega)}^2
            \\
            \ge &
       - 
        C \delta_{n-k}^2
   \|\varphi_{k+1}\|_{H^{1}(\Omega)}^2
			-
            \frac{\lambda}{4}
            \|\partial_{\boldsymbol{n}} \varphi_{k}\|_{L^{2}(\partial\Omega)}^2.
    \end{aligned}
\end{equation*}
Moreover, we apply the standard elliptic estimate of Neumann problem, see Proposition 2.10 of \cite{ern2004theory}:
\begin{equation}\label{standard applications of Lax-Milgram for the weak formulation,}
	\|\nabla\varphi_k\|_{L^2(\Omega)}^2
	\leq
	C(\|\Delta \varphi_k\|_{L^2(\Omega)}^2
 +
 \|\partial_{\boldsymbol{n}} \varphi_{k}\|_{L^{2}(\partial\Omega)}^2
 ).
\end{equation}
Note the fact $\delta_{n-k} \le \delta_{n-k-1} \le \delta$ for $0\le k\le n-2$, it follows that
\begin{equation*}
    \begin{aligned}
        2 \delta_{n-k} (\partial_{\boldsymbol{n}} \varphi_{k}, \varphi_{k+1})_{\partial\Omega}
             \ge &
        - 
        C \delta \big(\delta_{n-k-1}
   \|\Delta \varphi_{k+1}\|_{L^{2}(\Omega)}^2
   +
   \delta_{n-k}
   \| \varphi_{k+1}\|_{L^{2}(\Omega)}^2
   \big)\\
   & -
   C\delta \delta_{n-k-1}
   \|\partial_{\boldsymbol{n}} \varphi_{k+1}\|_{L^{2}(\partial\Omega)}^2
			-
            \frac{\lambda}{4}
            \|\partial_{\boldsymbol{n}} \varphi_{k}\|_{L^{2}(\partial\Omega)}^2.
    \end{aligned}
\end{equation*}
We apply the same argument to the last term in \eqref{boundary term for n case second order}:
\begin{equation*}
    \begin{aligned}
       -2 \delta_{n-k} ( \varphi_k, \partial_{\boldsymbol{n}} \varphi_{k+1})_{\partial\Omega}
            \ge &
        - 
        C \delta_{n-k}^2 
   \|\varphi_{k}\|_{H^{1}(\Omega)}^2
			-
            \frac{\lambda}{4}
            \|\partial_{\boldsymbol{n}} \varphi_{k+1}\|_{L^{2}(\partial\Omega)}^2\\
            \ge &
        - 
         C \delta \big(\delta_{n-k}
   \|\Delta \varphi_{k}\|_{L^{2}(\Omega)}^2
   +
   \delta_{n-k+1}
   \| \varphi_{k}\|_{L^{2}(\Omega)}^2
   \big)\\
   & -
   C \delta \delta_{n-k}
   \|\partial_{\boldsymbol{n}} \varphi_{k}\|_{L^{2}(\partial\Omega)}^2
			-
            \frac{\lambda}{4}
            \|\partial_{\boldsymbol{n}} \varphi_{k+1}\|_{L^{2}(\partial\Omega)}^2.
    \end{aligned}
\end{equation*}
Here, the inequality $\delta_{n-k}^2\leq\delta\delta_{n-k+1}$ for $0\le k\le n-2$ is used in the last line. Combining the estimates above and using notation \eqref{define I-second order}, we obtain
\begin{equation*}
\begin{aligned}
    2
  \sum_{k=0}^{n-2}\delta_{n-k}& \Big(
   (\partial_{\boldsymbol{n}} \varphi_{k}, \varphi_{k+1})_{\partial\Omega}
		-
   ( \varphi_{k},\partial_{\boldsymbol{n}} \varphi_{k+1})_{\partial\Omega}\Big)\\
  \ge& 
  - 
        C \delta 
       \mathcal{E}^*_{\delta,n}
       -
       \epsilon_{n}  \delta_{n+1}  \| \varphi_{0}\|_{L^2(\Omega)}^2
        -
        \big(\frac{\lambda}{2} +  C \delta \big)
        \sum_{k=0}^{n-1}
        \|\partial_{\boldsymbol{n}} \varphi_{k}\|_{L^{2}(\partial\Omega)}^2.
  \end{aligned}
\end{equation*}
Therefore the substituting into the boundary term \eqref{boundary term for n case second order} deduces the estimate \eqref{boundary estimate of D case n greater than 1-second order}.

\noindent\textbf{Case ii. Robin boundary}.
For the Robin boundary case, we utilize Young's inequality. The last two terms in \eqref{boundary term for n case second order} can be written as
\begin{equation}\label{rewrite boundary term Robin}
     \begin{aligned}
        &2 \delta_{n-k} (\partial_{\boldsymbol{n}} \varphi_{k}, \varphi_{k+1})_{\partial\Omega}
        -2 \delta_{n-k} ( \varphi_k, \partial_{\boldsymbol{n}} \varphi_{k+1})_{\partial\Omega}\\
        = &
        2 \delta_{n-k} ( \varphi_k + \partial_{\boldsymbol{n}} \varphi_{k}, \varphi_{k+1})_{\partial\Omega}
        -2 \delta_{n-k} ( \varphi_k, \varphi_{k+1} + \partial_{\boldsymbol{n}} \varphi_{k+1})_{\partial\Omega}.
    \end{aligned}
\end{equation}
The first term on the right-hand side of \eqref{rewrite boundary term Robin} yields the estimate
\begin{equation*}
    \begin{aligned}
        2 \delta_{n-k} (\varphi_k + \partial_{\boldsymbol{n}} \varphi_{k}, \varphi_{k+1})_{\partial\Omega}
        \ge &
        - 
        C \delta_{n-k}^2 
   \|\varphi_{k+1}\|_{L^{2}(\partial\Omega)}^2
			-
            \frac{\lambda}{4}
            \|\varphi_k + \partial_{\boldsymbol{n}} \varphi_{k}\|_{L^{2}(\partial\Omega)}^2
            \\
            \ge &
       - 
       C \delta_{n-k}^2 
   \|\varphi_{k+1}\|_{H^{1}(\Omega)}^2
			-
            \frac{\lambda}{4}
            \|\varphi_k + \partial_{\boldsymbol{n}} \varphi_{k}\|_{L^{2}(\partial\Omega)}^2.
    \end{aligned}
\end{equation*}
Note the fact $\delta_{n-k} \le \delta_{n-k-1} \le \delta$ for $0\le k\le n-2$, it follows that
\begin{equation*}
    \begin{aligned}
        2 \delta_{n-k} (\varphi_k + \partial_{\boldsymbol{n}} \varphi_{k}, \varphi_{k+1})_{\partial\Omega}
            \ge &
        - 
        C \delta \delta_{n-k-1}
   \|\varphi_{k+1}\|_{H^{1}(\Omega)}^2
			-
            \frac{\lambda}{4}
            \|\varphi_k + \partial_{\boldsymbol{n}} \varphi_{k}\|_{L^{2}(\partial\Omega)}^2.
    \end{aligned}
\end{equation*}
 Moreover, we apply the elliptic estimate of Robin problem:
\begin{equation}\label{elliptic estimate for Robin boundary}
	\|\varphi_k\|_{H^1(\Omega)}^2
	\leq
	C(\|\Delta \varphi_k\|_{L^2(\Omega)}^2
 +
 \|\varphi_k + \partial_{\boldsymbol{n}} \varphi_{k}\|_{L^{2}(\partial\Omega)}^2).
\end{equation}
Then it finally leads to
\begin{equation*}
    \begin{aligned}
       2 \delta_{n-k} (\varphi_k + \partial_{\boldsymbol{n}} \varphi_{k}, & \varphi_{k+1})_{\partial\Omega}
             \ge 
        - 
        C \delta \delta_{n-k-1}
   \|\Delta \varphi_{k+1}\|_{L^{2}(\Omega)}^2\\
   & -
   C \delta 
   \|\varphi_{k+1} + \partial_{\boldsymbol{n}} \varphi_{k+1}\|_{L^{2}(\partial\Omega)}^2
			-
            \frac{\lambda}{4}
            \|\varphi_k +\partial_{\boldsymbol{n}} \varphi_{k}\|_{L^{2}(\partial\Omega)}^2.
    \end{aligned}
\end{equation*}
We apply the same argument to the second term on the right-hand side of \eqref{rewrite boundary term Robin}:
\begin{equation*}
    \begin{aligned}
       & -2 \delta_{n-k} ( \varphi_k, \varphi_{k+1} +\partial_{\boldsymbol{n}} \varphi_{k+1})_{\partial\Omega}\\
            \ge &
        - 
         C \delta \delta_{n-k}
   \|\varphi_{k}\|_{H^{1}(\Omega)}^2
			-
            \frac{\lambda}{4}
            \|\varphi_{k+1} + \partial_{\boldsymbol{n}} \varphi_{k+1}\|_{L^{2}(\partial\Omega)}^2\\
            \ge &
        - 
        C\delta \delta_{n-k}
   \|\Delta \varphi_{k}\|_{L^{2}(\Omega)}^2
   -
   C \delta  
   \|\varphi_{k} + \partial_{\boldsymbol{n}} \varphi_{k}\|_{L^{2}(\partial\Omega)}^2
			-
            \frac{\lambda}{4}
            \|\varphi_{k+1} + \partial_{\boldsymbol{n}} \varphi_{k+1}\|_{L^{2}(\partial\Omega)}^2.
    \end{aligned}
\end{equation*}
Combining the estimates above, and using notation \eqref{define I-second order}, we obtain
\begin{equation*}
\begin{aligned}
   - 2
  \sum_{k=0}^{n-2}\delta_{n-k}& \Big(
   (\partial_{\boldsymbol{n}} \varphi_{k}, \varphi_{k+1})_{\partial\Omega}
		-
   ( \varphi_{k},\partial_{\boldsymbol{n}} \varphi_{k+1})_{\partial\Omega}\Big)\\
  \ge& 
  - 
       C \delta 
       \mathcal{E}^*_{\delta,n}
        -
        \big(\frac{\lambda}{2} +  C \delta \big)
        \sum_{k=0}^{n-1}
        \|\varphi_{k} + \partial_{\boldsymbol{n}} \varphi_{k}\|_{L^{2}(\partial\Omega)}^2.
  \end{aligned}
\end{equation*}

\noindent\textbf{Case iii. Dirichlet boundary}.
For the case of the Dirichlet boundary condition, the last two terms in \eqref{boundary term for n case} can be estimated by
\begin{equation*}
\begin{aligned}
     2 \delta_{n-k} (\partial_{\boldsymbol{n}} \varphi_{k}, \varphi_{k+1})_{\partial\Omega}
    \ge &
    -2 \delta
   \| \partial_{\boldsymbol{n}} \varphi_{k} \|_{L^{2}(\partial\Omega)}
   \| \varphi_{k+1} \|_{L^{2}(\partial\Omega)},\\
   -2 \delta_{n-k} ( \varphi_k, \partial_{\boldsymbol{n}} \varphi_{k+1})_{\partial\Omega}
        \ge &
    -2 \delta
   \| \partial_{\boldsymbol{n}} \varphi_{k+1} \|_{L^{2}(\partial\Omega)}
   \| \varphi_k\|_{L^{2}(\partial\Omega)}.
   \end{aligned}
\end{equation*}
Hence the estimate \eqref{estimate of boundary Dirichlet second order} follows.
\end{proof}

\subsection{Proof of Lemma \ref{Thm: second order}}
For the Neumann case, we substitute \eqref{boundary estimate of D case n greater than 1-second order} into \eqref{estimate of B_I n case second order} to derive that, for $\alpha = \mathrm{N}$,
  \begin{equation}\label{result of Neumann case second order}
		\begin{aligned}
			\mathcal{B}_\alpha(\boldsymbol{v};\boldsymbol{v})
			\ge&
   \big(1- \epsilon_{n} - C \delta \big) \mathcal{E}^*_{\delta,n}
   - 
  2 \epsilon_{n}  \delta_{n+1}  \| \varphi_{0}\|_{L^2(\Omega)}^2\\
     &   +
        \big(\frac{\lambda}{2} - C \delta \big)
        \sum_{k=0}^{n-1}
        \|\partial_{\boldsymbol{n}} \varphi_{k}\|_{L^{2}(\partial\Omega)}^2
        +
       \mu\Big|\int_\Omega \varphi_0\,\mathrm{d}\boldsymbol{x}\Big|^2.
		\end{aligned}
	\end{equation}
We utilize the inequality \eqref{standard applications of Lax-Milgram for the weak formulation,} and the Poincar\'e-Friedrichs inequality \eqref{Poincaré-Friedrichs inequality-2} to deduce
  \begin{equation*}
		\begin{aligned}
			2\epsilon_{n}  \delta_{n+1}  \| \varphi_{0}\|_{L^2(\Omega)}^2
			\le &
    \epsilon_{n}  \delta_{n+1}  C\Big(\|\Delta \varphi_0\|_{L^2(\Omega)}^2
 +
 \|\partial_{\boldsymbol{n}} \varphi_{0}\|_{L^{2}(\partial\Omega)}^2
 +
 \Big|\int_\Omega \varphi_0\,\mathrm{d}\boldsymbol{x}\Big|^2\Big)\\
   \le &
     C  \delta   \Big( \mathcal{E}^*_{\delta,n}
 +
 \|\partial_{\boldsymbol{n}} \varphi_{0}\|_{L^{2}(\partial\Omega)}^2
 +
 \Big|\int_\Omega \varphi_0\,\mathrm{d}\boldsymbol{x}\Big|^2\Big).
		\end{aligned}
	\end{equation*}
 Using the result above, it follows from \eqref{result of Neumann case second order} that
   \begin{equation*}
		\begin{aligned}
			\mathcal{B}_\alpha(\boldsymbol{v};\boldsymbol{v})
			\ge&
   \big(1- \epsilon_{n} -  C \delta \big) \mathcal{E}_{\delta,n}
   +
   \epsilon_{n}  \delta_{n+1}  \| \varphi_{0}\|_{L^2(\Omega)}^2\\
     &   +
        \big(\frac{\lambda}{2} -   C\delta\big)
        \sum_{k=0}^{n-1}
        \|\partial_{\boldsymbol{n}} \varphi_{k}\|_{L^{2}(\partial\Omega)}^2
        +
        (\mu- C \delta)\Big|\int_\Omega \varphi_0\,\mathrm{d}\boldsymbol{x}\Big|^2.
		\end{aligned}
	\end{equation*}
 Therefore \eqref{Neumann boundary condition-second order continuous} is derived for $\alpha = \mathrm{N}$. 

For the case of the Robin boundary, Combining \eqref{estimate of B_I n case second order} and \eqref{estimate of boundary Robin second order} to derive
\begin{equation*}
\begin{aligned}
\mathcal{B}_\alpha(\boldsymbol{v};\boldsymbol{v})
			\ge&
  \Big(1- \epsilon_{n}- C \delta \Big) \mathcal{E}^*_{\delta,n}\\
  & - 
   \epsilon_{n}  \delta_{n+1}  \| \varphi_{0}\|_{L^2(\Omega)}^2
       +
      \Big(\frac{\lambda}{2} - C \delta \Big)
       \sum_{k=0}^{n-1}
        \|\varphi_{k} + \partial_{\boldsymbol{n}} \varphi_k\|_{L^{2}(\partial\Omega)}^2.
        \end{aligned}
   \end{equation*}
   In addition. Applying the elliptic estimate \eqref{elliptic estimate for Robin boundary}, we have
   \begin{equation*}
   \begin{aligned}
       2\epsilon_{n}  \delta_{n+1}  \| \varphi_{0}\|_{L^2(\Omega)}^2
       \leq&
       \epsilon_{n}  \delta_{n+1}C(\|\Delta \varphi_0\|_{L^2(\Omega)}^2
 +
 \|\varphi_0 + \partial_{\boldsymbol{n}} \varphi_{0}\|_{L^{2}(\partial\Omega)}^2)\\
 \leq&
        C\delta (\mathcal{E}^*_{\delta,n}
 +
 \|\varphi_0 + \partial_{\boldsymbol{n}} \varphi_{0}\|_{L^{2}(\partial\Omega)}^2).
 \end{aligned}
   \end{equation*}
   Hence we conclude that
   \begin{equation*}
\begin{aligned}
\mathcal{B}_\alpha(\boldsymbol{v};\boldsymbol{v})
			\ge&
  \big(1- \epsilon_{n}-C \delta \big) \mathcal{E}^*_{\delta,n}\\
  & + 
   \epsilon_{n}  \delta_{n+1}  \| \varphi_{0}\|_{L^2(\Omega)}^2
       +
      \Big(\frac{\lambda}{2} - C \delta \Big)
       \sum_{k=0}^{n-1}
        \|\varphi_{k} + \partial_{\boldsymbol{n}} \varphi_k\|_{L^{2}(\partial\Omega)}^2.
        \end{aligned}
   \end{equation*}

For the case of Dirichlet boundary, we utilize \eqref{estimate of boundary Dirichlet second order}	and \eqref{estimate of B_I n case second order} to deduce that
 \begin{equation}
		\begin{aligned}
			\mathcal{B}_\alpha(\boldsymbol{v};\boldsymbol{v})
			\ge&
    (1- \epsilon_{n} ) \mathcal{E}^*_{\delta,n}
   - 
   \epsilon_{n}  \delta_{n+1}  \| \varphi_{0}\|_{L^2(\Omega)}^2\\
  & +
  \lambda \sum_{k=0}^{n-1} \Vert \varphi_{k}\Vert_{L^2(\partial\Omega)}^2
            -
            \sum_{k=0}^{n-1}
            \| \partial_{\boldsymbol{n}} \varphi_{k+1} \|_{L^{2}(\partial\Omega)}
            \| \varphi_k\|_{L^{2}(\partial\Omega)}.
		\end{aligned}
	\end{equation}
Applying Poincar\'e-Friedrichs inequality \eqref{Poincaré-Friedrichs inequality}, one can derive
    \begin{equation*}
    \begin{aligned}
        2\epsilon_{n}  \delta_{n+1}  \| \varphi_{0}\|_{L^2(\Omega)}^2
        \leq&
        \epsilon_{n}  \delta_{n+1} C(\|\nabla\varphi_{0}\|_{L^2(\Omega)}^2
        +
        \|\varphi_{0}\|_{L^2(\partial\Omega)}^2)\\
        \leq&
        C\delta(\|\Delta\varphi_{0}\|_{L^2(\Omega)}^2
        +
        \|\partial_n \varphi_{0} \|_{L^2(\partial\Omega)}\|\varphi_{0}\|_{L^2(\partial\Omega)}
        +
        \|\varphi_{0}\|_{L^2(\partial\Omega)}^2),
        \end{aligned}
    \end{equation*}
 where in the last line we used the elliptic estimate for Dirichlet boundary.
 Then by utilizing this fact $\| \varphi_k\|_{L^{2}(\partial\Omega)} 
     \le 
     \mathcal{B}^{\frac{1}{2}}_\alpha(\boldsymbol{v};\boldsymbol{v})$ , it follows that
 \begin{equation*}
     \begin{aligned}
         &\mathcal{B}_\alpha(\boldsymbol{v};\boldsymbol{v})
   +
   \sum_{k=0}^{n-1}
            \| \nabla\boldsymbol{\varphi}_{k} \|_{L^{2}(\partial\Omega)}
           \mathcal{B}^{\frac{1}{2}}_\alpha(\boldsymbol{v};\boldsymbol{v})\\
			\ge&
    (1- \epsilon_{n}-C\delta ) \mathcal{E}^*_{\delta,n}
   +
   \epsilon_{n}  \delta_{n+1}  \| \varphi_{0}\|_{L^2(\Omega)}^2
  +
  (\lambda-C\delta) \sum_{k=0}^{n-1} \Vert \varphi_{k}\Vert_{L^2(\partial\Omega)}^2.
     \end{aligned}
 \end{equation*}
 
 \section{Conclusion}
 This study delves into the error estimates of the Deep Mixed Residual method (MIM) in solving high-order elliptic equations with non-homogeneous boundary conditions, including Dirichlet, Neumann, and Robin conditions. Two types of loss functions of MIM, referred to as first-order and second-order least squares systems, are both considered.
 We apply a general approach in the notion of bilinear forms introduced in \cite{zeinhofer2023unified}, where one can utilize C\'ea's Lemma by performing boundedness and coercivity analysis. As a result, the total error is decomposed into three components: approximation error, generalization error, and optimization error. Through the Barron space theory and Rademacher complexity, an a priori error is derived regarding the training samples and network size that are exempt from the curse of dimensionality. Our results reveal that MIM significantly reduces the regularity requirements for activation functions compared to the deep Ritz method, implying the effectiveness of MIM in solving high-order equations.

 The main challenge in our analysis arises from the coercivity analysis of the low-order least squares systems (Equation \eqref{first order} and \eqref{second order}), specifically controlling the coupling terms. We employ a perturbation technique, selecting a special sequence of small parameters, to effectively bound these cross terms by decoupled terms. In addition, when considering the Dirichlet boundary condition, using an $L^2$ boundary penalty will lead to a loss of regularity of $3/2$. This implies that approximations in $H^2$ yield a posteriori estimate only in $H^{\frac{1}{2}}$. In this work, we utilize the idea in \cite{li2024priori} and derive a priori estimate in $H^1$ under the Dirichlet boundary condition. This is achieved by extending the C\'ea's Lemma and conducting a sup-linear coercivity analysis. 

 \section*{Acknowledgements}
  J. Chen is supported by the NSFC Major Research Plan - Interpretable and General-purpose Next-generation Artificial Intelligence (92370205) and NSFC 12425113. R. Du is supported by NSFC grant 12271360.
 The authors acknowledge partial support from the Austrian Science Fund (FWF), grant DOI 10.55776/P33010 and 10.55776/F65. This work has received funding from the European Research Council (ERC) under the European Union's Horizon 2020 research and innovation programme, ERC Advanced Grant NEUROMORPH, no. 101018153

\bibliographystyle{amsplain}
\bibliography{reference}

\providecommand{\bysame}{\leavevmode\hbox to3em{\hrulefill}\thinspace}
\providecommand{\MR}{\relax\ifhmode\unskip\space\fi MR }
\providecommand{\MRhref}[2]{%
  \href{http://www.ams.org/mathscinet-getitem?mr=#1}{#2}
}
\providecommand{\href}[2]{#2}
\begin{thebibliography}{10}

\bibitem{andersson1998solution}
Lars-Erik Andersson, Tommy Elfving, and Gene~H Golub, \emph{Solution of
  biharmonic equations with application to radar imaging}, Journal of
  Computational and Applied Mathematics \textbf{94} (1998), no.~2, 153--180.

\bibitem{barron1993universal}
Andrew~R Barron, \emph{Universal approximation bounds for superpositions of a
  sigmoidal function}, IEEE Transactions on Information theory \textbf{39}
  (1993), no.~3, 930--945.

\bibitem{cai2020deep}
Zhiqiang Cai, Jingshuang Chen, Min Liu, and Xinyu Liu, \emph{Deep least-squares
  methods: An unsupervised learning-based numerical method for solving elliptic
  {PDEs}}, Journal of Computational Physics \textbf{420} (2020), 109707.

\bibitem{cai1994first}
Zhiqiang Cai, R~Lazarov, Thomas~A Manteuffel, and Stephen~F McCormick,
  \emph{First-order system least squares for second-order partial differential
  equations: Part {I}}, SIAM Journal on Numerical Analysis \textbf{31} (1994),
  no.~6, 1785--1799.

\bibitem{duan2007exact}
WH~Duan and Chien~Ming Wang, \emph{Exact solutions for axisymmetric bending of
  micro/nanoscale circular plates based on nonlocal plate theory},
  Nanotechnology \textbf{18} (2007), no.~38, 385704.

\bibitem{yu2018deep}
Weinan E and Bing Yu, \emph{The deep ritz method: a deep learning-based
  numerical algorithm for solving variational problems}, Communications in
  Mathematics and Statistics \textbf{6} (2018), no.~1, 1--12.

\bibitem{ern2004theory}
Alexandre Ern and Jean-Luc Guermond, \emph{Theory and practice of finite
  elements}, vol. 159, Springer, 2004.

\bibitem{girault2012finite}
Vivette Girault and Pierre-Arnaud Raviart, \emph{Finite element methods for
  {Navier--Stokes} equations: theory and algorithms}, vol.~5, Springer Science
  \& Business Media, 2012.

\bibitem{han2018solving}
Jiequn Han, Arnulf Jentzen, and Weinan E, \emph{Solving high-dimensional
  partial differential equations using deep learning}, Proceedings of the
  National Academy of Sciences \textbf{115} (2018), no.~34, 8505--8510.

\bibitem{hong2021priori}
Qingguo Hong, Jonathan~W Siegel, and Jinchao Xu, \emph{A priori analysis of
  stable neural network solutions to numerical {PDEs}}, arXiv preprint
  arXiv:2104.02903 (2021).

\bibitem{khain2008generalized}
Evgeniy Khain and Leonard~M Sander, \emph{Generalized {Cahn-Hilliard} equation
  for biological applications}, Physical Review E—Statistical, Nonlinear, and
  Soft Matter Physics \textbf{77} (2008), no.~5, 051129.

\bibitem{landau2012theory}
Lev~Davidovich Landau, LP~Pitaevskii, Arnol'd~Markovich Kosevich, and
  Evgenii~Mikhailovich Lifshitz, \emph{Theory of elasticity: volume 7}, vol.~7,
  Elsevier, 2012.

\bibitem{li2024priori}
Lingfeng Li, Xue-Cheng Tai, Jiang Yang, and Quanhui Zhu, \emph{A priori error
  estimate of deep mixed residual method for elliptic {PDEs}}, Journal of
  Scientific Computing \textbf{98} (2024), no.~2, 44.

\bibitem{li2024two}
Yuanyuan Li, Shuai Lu, Peter Math{\'e}, and Sergei~V Pereverzev,
  \emph{Two-layer networks with the {$ ReLU^k$} activation function: Barron
  spaces and derivative approximation}, Numerische Mathematik \textbf{156}
  (2024), no.~1, 319--344.

\bibitem{lions2012non}
Jacques~Louis Lions and Enrico Magenes, \emph{Non-homogeneous boundary value
  problems and applications: Vol. 1}, vol. 181, Springer Science \& Business
  Media, 2012.

\bibitem{lu2021priori}
Yulong Lu, Jianfeng Lu, and Min Wang, \emph{A priori generalization analysis of
  the deep ritz method for solving high dimensional elliptic partial
  differential equations}, Conference on learning theory (2021), 3196--3241.

\bibitem{lyu2020enforcing}
Liyao Lyu, Keke Wu, Rui Du, and Jingrun Chen, \emph{Enforcing exact boundary
  and initial conditions in the deep mixed residual method}, CSIAM Transactions
  on Applied Mathematics \textbf{2} (2021), no.~4, 748--775.

\bibitem{lyu2022mim}
Liyao Lyu, Zhen Zhang, Minxin Chen, and Jingrun Chen, \emph{M{IM}: A deep mixed
  residual method for solving high-order partial differential equations},
  Journal of Computational Physics \textbf{452} (2022), 110930.

\bibitem{mishra2021estimates}
Siddhartha Mishra and Roberto Molinaro, \emph{Estimates on the generalization
  error of physics informed neural networks ({PINNs}) for approximating a class
  of inverse problems for {PDEs}}, IMA Journal of Numerical Analysis
  \textbf{42} (2022), no.~2, 981--1022.

\bibitem{mishra2023estimates}
\bysame, \emph{Estimates on the generalization error of physics-informed neural
  networks for approximating {PDEs}}, IMA Journal of Numerical Analysis
  \textbf{43} (2023), no.~1, 1--43.

\bibitem{muller2021error}
Johannes M{\"u}ller and Marius Zeinhofer, \emph{Error estimates for the
  variational training of neural networks with boundary penalty}, arXiv
  preprint arXiv:2103.01007 (2021).

\bibitem{muller2022error}
\bysame, \emph{Error estimates for the deep ritz method with boundary penalty},
  Mathematical and Scientific Machine Learning (2022), 215--230.

\bibitem{pmlr-v190-muller22b}
\bysame, \emph{Notes on exact boundary values in residual minimisation},
  Mathematical and Scientific Machine Learning (2022), 231--240.

\bibitem{nevcas1967methodes}
Jind{\v{r}}ich Ne{\v{c}}as, \emph{Les m{\'e}thodes directes en th{\'e}orie des
  {\'e}quations elliptiques}, Academia (1967).

\bibitem{raissi2019physics}
Maziar Raissi, Paris Perdikaris, and George~E Karniadakis,
  \emph{Physics-informed neural networks: A deep learning framework for solving
  forward and inverse problems involving nonlinear partial differential
  equations}, Journal of Computational physics \textbf{378} (2019), 686--707.

\bibitem{schechter1963p}
Martin Schechter, \emph{On {$L^p$} estimates and regularity {II}}, Mathematica
  Scandinavica \textbf{13} (1963), no.~1, 47--69.

\bibitem{siegel2023greedy}
Jonathan~W Siegel, Qingguo Hong, Xianlin Jin, Wenrui Hao, and Jinchao Xu,
  \emph{Greedy training algorithms for neural networks and applications to
  {PDEs}}, Journal of Computational Physics \textbf{484} (2023), 112084.

\bibitem{sirignano2018dgm}
Justin Sirignano and Konstantinos Spiliopoulos, \emph{D{GM}: A deep learning
  algorithm for solving partial differential equations}, Journal of
  computational physics \textbf{375} (2018), 1339--1364.

\bibitem{vahab2022physics}
Mohammad Vahab, Ehsan Haghighat, Maryam Khaleghi, and Nasser Khalili, \emph{A
  physics-informed neural network approach to solution and identification of
  biharmonic equations of elasticity}, Journal of Engineering Mechanics
  \textbf{148} (2022), no.~2, 04021154.

\bibitem{weinan2019barron}
E~Weinan, Chao Ma, and Lei Wu, \emph{The {B}arron space and the flow-induced
  function spaces for neural network models}, Constructive Approximation.
  \textbf{55} (2022), no.~1, 369--406.

\bibitem{zang2020weak}
Yaohua Zang, Gang Bao, Xiaojing Ye, and Haomin Zhou, \emph{Weak adversarial
  networks for high-dimensional partial differential equations}, Journal of
  Computational Physics \textbf{411} (2020), 109409.

\bibitem{zeinhofer2023unified}
Marius Zeinhofer, Rami Masri, and Kent-Andr{\'e} Mardal, \emph{A unified
  framework for the error analysis of physics-informed neural networks}, arXiv
  preprint arXiv:2311.00529 (2023).

\end{thebibliography}

\appendix
\section{Neural-Network Approximation}\label{Sec:Neural-Network Approximation}
\subsection{The proof of Lemma \ref{Lem: approximation of ReCU network}}
For simplicity, we will only provide the proof for the case when the output function is one-dimensional here. The higher-dimensional case can be directly generalized.
From \cite{lu2021priori},
we recall for any $s\in\mathbb{N}$, the spectral Barron space $\mathcal{B}^s(\Omega)$ given by
\begin{equation*}
    \mathcal{B}^s(\Omega)
    :=
    \Big\{
    u\in L^1(\Omega)\,\Big|\,
    \sum_{k\in \mathbb{N}_0^d}(1+\pi^s|k|_1^s)|\hat{u}(k)|
    <\infty
    \Big\},
\end{equation*}
with associated norm $\|u\|_{\mathcal{B}^s(\Omega)}:=\sum_{k\in \mathbb{N}_0^d}(1+\pi^s|k|_1^s)|\hat{u}(k)|$, where $\{\hat{u}(k)\}_{k\in \mathbb{N}_0^d}$ is the expansion cofficients of $u$ under the cosine basis functions. And for any function $u\in H^1(\Omega)$ admits the expansion 
\begin{equation}
    u(\boldsymbol{x})
    :=
    \hat{u}(0)
    +
    \int
    g(\boldsymbol{x},k)\mu(\mathrm{d}k),
\end{equation}
where $\mu(\mathrm{d}k)$ is the probability measure on $\mathbb{N}_0^d/ \{\boldsymbol{0}\}$ define by
\begin{equation*}
    \mu(\mathrm{d}k)
    =
    \sum_{\mathbb{N}_0^d/ \{\boldsymbol{0}\}} \frac{1}{Z_u}|\hat{u}(k)|
    (1+\pi^4|k|_1^4)\delta(\mathrm{d}k),
\end{equation*}
with normalizing constant $Z_u=\sum_{\mathbb{N}_0^d/ \{\boldsymbol{0}\}} |\hat{u}(k)|(1+\pi^4|k|_1^4)\leq \|u\|_{\mathcal{B}^4(\Omega)}$, and 
\begin{equation*}
   g(\boldsymbol{x},k)
   =
   \frac{Z_u}{1+\pi^4|k|_1^4}\cdot\frac{1}{2^d}\sum_{\xi\in\Xi}\cos\big(\pi(\boldsymbol{k}_\xi\cdot\boldsymbol{x}+\theta_k)\big),
\end{equation*}
with $\theta(k)\in \{0,1\}$ and $\boldsymbol{k}_\xi=(k_1\xi_1,\cdots k_d\xi_d)$.
Let us define the function class:
\begin{equation*}
    \mathcal{F}_{\cos}(B)
    :=
    \left\{
    \frac{\gamma}{1+\pi^4|k|_1^4}\cos\big(\pi(\boldsymbol{k}\cdot\boldsymbol{x})+b\big)\,\Big|\, k\in \mathbb{Z}^d/\{\boldsymbol{0}\}, |\gamma|\leq B, b\in \{0,1\}
    \right\},
\end{equation*}
where  $B>0$ is a constant. 
Hence if $u\in \mathcal{B}^4(\Omega)$, then $\Bar{u}:=u-\hat{u}(0)$ lies in the $H^1$-closure of the convex hull of $\mathcal{F}_{\cos}(B)$ with $B=\|u\|_{\mathcal{B}^4(\Omega)}$. Since the $H^1$-norm of any function in $\mathcal{F}_{\cos}(B)$ is bounded by $B$, apply the following conclusion:
\begin{lemma}\cite{barron1993universal}
  Suppose $u$ belongs to the closure of the convex hull of a set $\mathcal{G}$ in a Hilbert space. Let the Hilbert norm of each element of $\mathcal{G}$ be bounded by $B>0$. Then for every $m\in \mathbb{N}$, there exists $\{g_i\}_{i=1}^m\subset\mathcal{G}$ and $\{c_i\}_{i=1}^m\subset[0,1]$ with $\sum_{i=1}^m c_i=1$ such that
  \begin{equation*}
      \Big\|
      u-\sum_{i=1}^m c_ig_i
      \Big\|^2
      \leq
      \frac{B^2}{m}.
  \end{equation*}
\end{lemma}
 Therefore it yields the following theorem:
 \begin{theorem}\label{Thm: approximation apply to Barron space}
     Let $u\in \mathcal{B}^4$. Then there exists $\phi_m$ which is a convex combination of $m$ functions in $\mathcal{F}_{\cos}(B)$ with $B=\|u\|_{\mathcal{B}^4(\Omega)}$ such that
     \begin{equation*}
         \|u-\hat{u}(0)-\phi_m\|_{H^1(\Omega)}^2
         \leq
         \frac{\|u\|_{\mathcal{B}^4(\Omega)}^2}{m}.
     \end{equation*}
 \end{theorem}

Next, we will give the reduction to the ReCU function. Notice that every function in $\mathcal{F}_{\cos}(B)$ is the composition of the one dimensional function $g$ defined on $[-1,1]$ by
\begin{equation}
    g(z)
    =
    \frac{\gamma}{1+\pi^4|k|_1^4}\cos\big(\pi(|k|_1z+b)\big),
\end{equation}
with $k\in \mathbb{Z}^d/\{\boldsymbol{0}\}$, $|\gamma|\leq B$ and $b\in \{0,1\}$, and a linear function $z=\boldsymbol{w}\cdot\boldsymbol{x}$ with $\boldsymbol{w}=\boldsymbol{k}/|\boldsymbol{k}|_1$. It is clear that $g\in C^4([-1,1]) $ satisfies 
\begin{equation*}
    \|g^{(s)}\|_{L^\infty[-1,1]}
    \leq
    |\gamma|
    \leq
    B
    \quad
    \mbox{for}\ s=0,1,\cdots,4.
\end{equation*}
\begin{lemma}
    Let $g\in C^4([-1,1])$ with $\|g^{(s)}\|_{L^\infty([-1,1])}\leq B$ for $s=0,1,\cdots,4$. Assume $g^{(s)}=0$ for $s=1,2,3$. Let $\{z_j\}_{j=0}^{2m}$ be a partition of $[-1,1]$ with $z_0=-1,z_m=0,z_{2m}=1$ and $z_{j+1}-z_j=h=1/m$ for each $j=0,\cdots,2m-1$. Then there exists a two-layer ${\rm ReCU}$ networks $g_m$ of the form 
    \begin{equation}
        g_m(z)
        =
        c
        +
        \sum_{i=1}^{2m+5}a_i {\rm ReCU}(\epsilon_i z-b_i),
    \end{equation}
    with $c=g(0),\ b_i\in [-1,1],\ \epsilon_i\in \{-1,1\}\ \mbox{for}\ i=1,\cdots,2m+5$, and  $\sum_{i=1}^{2m+5}|a_i|\leq 8B$ such that
    \begin{equation*}
        \|g-\hat{g}_m\|_{H^1(\Omega)}
        \leq
        \frac{6B}{\sqrt{m}}.
    \end{equation*}
\end{lemma}
\begin{proof}
    From \cite{li2024priori}, There exists a two-layer ReQU activated neural networks $\hat{g}_m(z)$ of the form
    \begin{equation}\label{ReQU activated neural netwoeks}
        \hat{g}_m(z)
        =
        c
        +
        \sum_{i=0}^{2m+3}\hat{a}_i \mbox{ReQU}(\epsilon_i z-\hat{b}_i),
    \end{equation}
    with $c=g(0), \hat{b}_i\in [-1,1]$, $\epsilon\in \{-1,1\},\ \mbox{for}\  i=0,\cdots,2m+3$, and $\sum_{i=0}^{2m+3}|\hat{a}_i|\leq 8B$ such that
    \begin{equation*}
        \|g-\hat{g}_m\|_{H^1(\Omega)}
        \leq
        \frac{5B}{\sqrt{m}}.
    \end{equation*}
    In addition, an example of coefficients $\{\hat{a}_i\}_{i=0}^{2m+3}$ are given as
    \begin{equation*}
    \hat{a}_i=
    \left\{
        \begin{aligned}
            &\frac{\tilde{a}_{i+1}}{4h}&i=0,1,\\
            &\frac{\tilde{a}_{i+1}-\tilde{a}_{i-1}}{4h}& 2\leq i\leq m-1,\\
            &-\frac{\tilde{a}_{i-1}}{4h}& i=m,m+1,\\
            &\frac{\tilde{a}_{i-1}}{4h}& i=m+2,m+3,\\
            &\frac{\tilde{a}_{i-1}-\tilde{a}_{i-3}}{4h}& m+4\leq i\leq 2m+1,\\
            &-\frac{\tilde{a}_{i-3}}{4h}&i=2m+2,2m+3,
        \end{aligned}
         \right.
    \end{equation*}
    with
    \begin{equation*}
   \tilde{a}_i=
    \left\{
        \begin{aligned}
          &\frac{g(z_{m+1})-g(z_{m})}{h}, &i=m+1,\\
          &\frac{g(z_{m-1})-g(z_{m})}{h}, &i=m,\\
          &\frac{g(z_{i})-2g(z_{i-1})+g(z_{i-2})}{h}, &i>m+1,\\
          &\frac{g(z_{i-1})-2g(z_{i})+g(z_{i+1})}{h}, &i<m.
        \end{aligned}
        \right.
    \end{equation*}
   And setting $\epsilon_i=-1, \hat{b}_i=-z_{i}$ for $i=0,\cdots,m+1$ and $\epsilon_i=1, \hat{b}_i=z_{i-3}$ for $i=m+2,\cdots,2m+3$.
    Notice that
    \begin{equation*}
        6\delta^2 \mbox{ReLU}(z)+6\delta \mbox{ReQU}(z)
        =
      \mbox{ReCU}(z+\delta)
        -
      \mbox{ReCU}(z-\delta)
        +
       e(z,\delta)
    \end{equation*}
    with
    \begin{equation*}
    \hat{e}(z;\delta)
    =
        \left\{
        \begin{aligned}
            &-(z+\delta)^2,\quad &-\delta<z\leq 0,\\
            &-(z+\delta)^2,\quad &0<z\leq \delta,\\
            &0,\quad&\mbox{otherwise},
        \end{aligned}
        \right.
    \end{equation*}
    and
    \begin{equation*}
    e(z;\delta)
    =
        \left\{
        \begin{aligned}
            &-(z+\delta)^3,\quad &-\delta<z\leq 0,\\
            &-(z+\delta)^3,\quad &0<z\leq \delta,\\
            &0,\quad&\mbox{otherwise}.
        \end{aligned}
        \right.
    \end{equation*}
Hence we have
\begin{equation}
    \begin{aligned}
        6\delta \mbox{ReQU}(z)
        =&
      \mbox{ReCU}(z+\delta)
        -
      \mbox{ReCU}(z-\delta)\\
      &-
      6\delta^2 \mbox{ReLU}(z)
        +
       e(z,\delta)-\frac{3}{2}\delta\hat{e}(z;\delta).
    \end{aligned}
\end{equation}
For fixed $\delta=h$, we obtain a ReCU activated neural networks $g_m(z)$ by approximating the ReQU active function in \eqref{ReQU activated neural netwoeks}:
\begin{equation}\label{ReCU activated neural network}
    g_m(z)
    =
    g(0)
    +
    \sum_{i=0}^{2m+3}\frac{\hat{a}_i}{6h}
    \big(\mbox{ReCU}(\epsilon_i z-\hat{b}_i+h)
        -
      \mbox{ReCU}(\epsilon_i z-\hat{b}_i-h)\big).
\end{equation}
It is easy to deduce that directly
\begin{equation}
    \begin{aligned}
        \|g-g_m\|_{H^1(\Omega)}
        \leq&
        \|g-\hat{g}_m\|_{H^1(\Omega)}
        +
        \|\hat{g}_m-g_m\|_{H^1(\Omega)}\\
        \leq&
        \frac{5B}{\sqrt{m}}+\frac{B}{\sqrt{m}}=\frac{6B}{\sqrt{m}}.
    \end{aligned}
\end{equation}
Rearrange \eqref{ReCU activated neural network} in the form of 
\begin{equation*}
  g_m(z)
  =
  g(0)
  +
  \sum_{i=1}^{m+2}a_i \mbox{ReCU}(z_i-z)
  +
  \sum_{i=m+3}^{2m+5}a_i \mbox{ReCU}(z-z_{i-5}),
\end{equation*}
where 
\begin{equation*}
a_i=
\left\{
    \begin{aligned}
      & \frac{\hat{a}_{i+1}-\hat{a}_{i-1}}{9h}, &1\leq i\leq m,\\
      & -\frac{\hat{a}_{i-1}}{9h},& i=m+1,m+2,\\
      & \frac{\hat{a}_{i-1}}{9h}, &i=m+3,m+4,\\
      & \frac{\hat{a}_{i-1}-\hat{a}_{i-3}}{9h},& m+5\leq i\leq 2m+4,\\
      & -\frac{\hat{a}_{i-3}+\hat{a}_{i-2}}{9h},& i=2m+5.
    \end{aligned}
    \right.
\end{equation*}
Furthermore, by the mean value theorem, $a_{m+1}=-\frac{\hat{a}_{m}}{9h}\leq |g^{(4)}h|\leq Bh$, which holds because of the assumption $g^{(s)}=0$ for $s=1,2,3$. It is similar to estimating the other coefficients.

Finally, by setting $\epsilon_i=-1, b_i=-z_i$ for $i=1,\cdots,m+2$ and $\epsilon_i=1, b_i=z_{i-5}$ for $i=m+3,\cdots,2m+5$, it completes the proof of the lemma. 
\end{proof}
Hence we have the following proposition:
\begin{proposition}\label{Pro:F ReCU B}
    Define the function class
    \begin{equation*}
        \mathcal{F}_{ \rm ReCU}(B)
        :=
        \left\{
        c+\gamma {\rm ReCU}(\boldsymbol{w}\cdot\boldsymbol{x}-t),
        |c|\leq 2B, |\boldsymbol{w}|_1=1, |t|\leq 1, |\gamma|\leq 4B
        \right\}.
    \end{equation*}
    Then for any constant $\tilde{c}$ such that $|\tilde{c}|\leq B$, the set $\tilde{c}+\mathcal{F}_{\cos}(B)$ is in the $H^1$-closure of the convex hull of $\mathcal{F}_{ \rm ReCU}(B)$.
\end{proposition}

With Proposition \ref{Pro:F ReCU B}, we are ready to give the proof of Lemma \ref{Lem: approximation of ReCU network}.
\begin{proof}[Proof of Lemma \ref{Lem: approximation of ReCU network}]
    Observe that if $u\in \mathcal{F}_{ \rm ReCU}(B)$, then
    \begin{equation*}
        \|u\|_{H^1(\Omega)}^2
        \leq
        (c+8\gamma)^2
        +
        144\gamma^2.
    \end{equation*}
    Therefore Lemma \ref{Lem: approximation of ReCU network} follows directly from Theorem \ref{Thm: approximation apply to Barron space} and Proposition \ref{Pro:F ReCU B}.
\end{proof}

\subsection{The proof of lemma \ref{Lem: Rademacher complexity for second-order system}}\label{Subsec: The proof of generalization error}
Next, we will effectively use Rademacher complexity to derive an upper bound on the generalization error. First, let us review some important properties of Rademacher complexity from \cite{li2024priori}:
\begin{lemma}\label{Lem: properties of Rademacher complexity}
    Let $\mathcal{F}, \mathcal{G}$ be function classes and $a,b$ be constants. Then
    \begin{enumerate}[(1)]
        \item\label{property 1}$R_{N}(\mathcal{F}+\mathcal{G})\leq R_{N}(\mathcal{F})+R_{N}(\mathcal{G}) $.
        \item\label{property 2} $R_{N}(a\mathcal{F})=|a|R_{N}(\mathcal{F})$.
        \item\label{property 3} Assume $g$ is a fixed function and $\|g\|_{L^\infty}\leq b$, then $R_\mathrm{N}(\mathcal{G})\leq \frac{b}{\sqrt{N}}$.
        \item\label{property 4} Assume that $\sigma:\mathbb{R}\mapsto\mathbb{R}$ is $l$-Lipschitz with $\sigma(0)=0$, then $R_{N}(\sigma(\mathcal{F}))\leq 2l R_{N}(\mathcal{F})$.
        \item\label{property 5} $R_{N}(\mathcal{F}^2)\leq 4\sup_{f\in \mathcal{F}}\|f\|_{L^\infty}R_{N}(\mathcal{F})$.
        \item\label{property 6} $R_{N}(\mathcal{F}\mathcal{G})\leq 6\sup_{f\in \mathcal{F}\cup\mathcal{G}} \|f\|_{L^\infty}(R_{N}(\mathcal{F})+R_{N}(\mathcal{G}))$.
    \end{enumerate}
\end{lemma}
With the calculation rules prepared, we are ready to estimate the complexity of the neural network function classes.
\begin{lemma} \label{Lem: linear transformation function class}
\cite{li2024priori}
    Let $\mathcal{G}$ be the linear transformation function class defined by 
    \begin{equation}
        \mathcal{G}
        :=
        \{\boldsymbol{\omega}\cdot\boldsymbol{x}+b
        \,|\,
        \|\boldsymbol{w}\|_2=1,\,|b|\leq 1
        \}.
    \end{equation}
    Then we have
    \begin{equation}
        R_N(\mathcal{G})
        \leq 
        \frac{\sqrt{2d\log d}+1}{\sqrt{N}}.
    \end{equation}
\end{lemma}
Similar to the function class $\mathcal{F}_{{\rm ReCU},m}$ in \cite{li2024priori}, we give an estimate of the Rademacher complexity for the function class $\mathcal{F}_{{\rm ReCU},m}$. 
The estimate of the complexity of a two-layer neural network depends on the activation function. Since $\sigma=\rm ReCU$ is three times three times continuously differentiable, we can make the following assumptions:
\begin{equation*}
    \sup |\sigma^{(k)}|\leq \ell_k,
    \quad
    k=0,1,2,3.
\end{equation*}
\begin{lemma}\label{Lemm: The Rademacher complexity for ReCU}
The Rademacher complexity of $\mathcal{F}_{{\rm ReCU},m}$ is bounded by
     \begin{equation}
        R_{N}(\mathcal{F}_{{\rm ReCU},m})
        \leq
        \frac{C\|u\|_{\mathcal{B}^{2n+3}}}{\sqrt{N}},
    \end{equation}
    where $C>0$ is a constant dependent on $d$.
\end{lemma}
\begin{proof}
    Using properties \eqref{property 1} and \eqref{property 3}, the Rademacher complexity of $\mathcal{F}_{{\rm ReCU},m}$ is broken down into the sum of the Rademacher complexity of each neuron, i.e. 
    \begin{equation*}
        R_N(\mathcal{F}_{{\rm ReCU},m})
        \leq
        \frac{2\|u\|_{\mathcal{B}^{2n+3}}}{\sqrt{N}}
        +
        \sum_{i=1}^m |\boldsymbol{a}_i|R_N\big(\sigma(\mathcal{G})\big).
    \end{equation*}
    Since $\sigma$ is $\ell_1$-Lipschitz and $\sigma(0)=0$, property \eqref{property 4} tells us
    \begin{equation*}
        R_N(\sigma\big(\mathcal{G})\big)
        \leq
        2\ell_1 R_N(\mathcal{G}).
    \end{equation*}
    Hence combining lemma \ref{Lem: linear transformation function class}, we can conclude
    \begin{equation*}
         R_N(\mathcal{F}_{{\rm ReCU},m})
        \leq
        \frac{(2+8\ell_1+8\ell_1\sqrt{2d\log d})\|u\|_{\mathcal{B}^{2n+3}}}{\sqrt{N}}.
    \end{equation*}
 \end{proof}
 With the preparation above, we will now prove Lemma \ref{Lem: Rademacher complexity for second-order system}.
 \begin{proof}[The proof of lemma \ref{Lem: Rademacher complexity for second-order system}]
By the definition of $\mathcal{F}_{{\rm ReCU},m}$ neural network, for $1\leq j\leq d$, we deduce that the first-order derivation is 
\begin{equation*}
    \frac{\partial \boldsymbol{v}_\theta}{\partial x_j}
    =
     \sum_{i=1}^{m}\boldsymbol{a}_i\cdot W_{i,j}\sigma^{(1)}(W_i\boldsymbol{x}+\boldsymbol{b}_i),
\end{equation*}
and second-order derivation is
\begin{equation*}
    \frac{\partial^2 \boldsymbol{v}_\theta}{\partial^2 x_j}
    =
     \sum_{i=1}^{m}\boldsymbol{a}_i\cdot |W_{i,j}|^2\sigma^{(2)}(W_i\boldsymbol{x}+\boldsymbol{b}_i).
\end{equation*}
Hence we have
\begin{equation}\label{V theta L infty estimate}
    \Big\|\frac{\partial^k \boldsymbol{v}_\theta}{\partial^k x_j}\Big\|_{L^\infty(\Omega)}
    \leq
    C \|u\|_{\mathcal{B}^{2n+3}},\quad \mbox{for}\  k=0,1,2.
\end{equation}
Recall the definition \eqref{Def: second-order matrix} of $\mathcal{P}^*$, it follows that
    \begin{equation}
        |\mathcal{P}^* \boldsymbol{v}_\theta-\boldsymbol{f}|^2
        =
        |\Delta \varphi_{n-1}-f|^2
        +
        \sum_{i=0}^{n-2}|\varphi_{i+1}-\Delta \varphi_{i}|^2.
     \end{equation}
 It follows that from property \eqref{property 5} 
 \begin{equation}
 \begin{aligned}
     R_N(L^*)
     \leq&
     4\sup_{\boldsymbol{v}_\theta\in \mathcal{F}_{{\rm ReCU},m}}\|\mathcal{P}^* \boldsymbol{v}_\theta-\boldsymbol{f}\|_{L^\infty(\Omega)}
     \Big(R_N(f)\\
         &+
         \sum_{i=1}^m \sum_{j=1}^d|\boldsymbol{a}_i| |W_{i,j}|^2 R_N\big(\sigma^{(2)}(\mathcal{G})\big)
         \Big).
         \end{aligned}
 \end{equation}
 According to the definition of the neural network and \eqref{V theta L infty estimate}, we have
    \begin{equation*}
    \sup_{\boldsymbol{v}_\theta\in \mathcal{F}_{{\rm ReCU},m}}
    \|\mathcal{P}^* \boldsymbol{v}_\theta-\boldsymbol{f}\|_{L^\infty(\Omega)}
  \leq
   C\|u^*\|_{\mathcal{B}^{2m+3}},
    \end{equation*}
    and
    \begin{equation*}
    \begin{aligned}
        \sum_{j=1}^d\sum_{i=1}^m |\boldsymbol{a}_i| |W_{i,j}|^2R_{N}(\sigma^{(2)}(\mathcal{G}))
        \leq&
        \sum_{i=1}^m d|a_i| \|W_{i}\|_2^2 R_{N}(\sigma^{(2)}(\mathcal{G}))\\
        \leq&
       Cd\|u^*\|_{\mathcal{B}^{2m+3}}R_{N}(\mathcal{G}),
        \end{aligned}
    \end{equation*}
    where the last inequality holds by the property \eqref{property 4}. Combining lemma \ref{Lem: linear transformation function class}, it yields that
    \begin{equation}
         R_N(L^*)
         \leq
         C\sqrt{\frac{2d^3\log d}{N}}\|u^*\|_{\mathcal{B}^{2m+3}}^2.
    \end{equation}
    Moreover, we estimate the function classes $L^*_{b,\alpha}$ with respect to the boundary conditions. In the same manner, for $\alpha=D$, we have 
    \begin{equation}\label{b-lemma rademacher}
   R_{\widehat{N}}(L_{\alpha}^*) 
   \leq
   4\sup_{\boldsymbol{v}_\theta\in \mathcal{F}_{{\rm ReCU},m}}
   \|S_\alpha \boldsymbol{v}_\theta-\boldsymbol{g}_\alpha\|_{L^\infty}
   \big( R_{\widehat{N}}(\mathcal{F}_{{\rm ReCU},m})
   +
   R_{\widehat{N}}(\boldsymbol{g}_\alpha) \big),
\end{equation}
where use the fact \eqref{V theta L infty estimate} to get
\begin{equation}
    \sup_{\boldsymbol{v}_\theta\in \mathcal{F}_{{\rm ReCU},m}}
   \|S_\alpha \boldsymbol{v}_\theta-\boldsymbol{g}_\alpha\|_{L^\infty}
   \leq
   C\|u^*\|_{\mathcal{B}^{2m+3}}.
\end{equation}
Therefore, Combining lemma \ref{Lem: linear transformation function class} and \ref{Lemm: The Rademacher complexity for ReCU}, we deduce that
\begin{equation}
     R_{\widehat{N}}(L_{\alpha}^*) 
     \leq
     C\sqrt{\frac{2d\log d}{\widehat{N}}}\|u^*\|_{\mathcal{B}^{2m+3}}^2
\end{equation}
Hence  taking $\widehat{N}=\frac{N}{d^2}$ will not
change the upper bound of the Rademacher complexity. The situation for other boundary conditions can be obtained similarly. 
    \end{proof}

\end{document}